\documentclass[a4paper,12pt]{scrartcl}

\usepackage[utf8]{inputenc}

\usepackage{amsthm,amsmath,amsfonts,amssymb}
\numberwithin{equation}{section}  
\usepackage{mathrsfs}
\usepackage{todonotes}

\usepackage{authblk}

\usepackage{tikz-cd}

\usepackage{mathtools}	 
\usepackage{bm}	  
\usepackage{bbm}	 

\usepackage{tensor}	 
\usepackage[colorlinks=true]{hyperref} 
\usepackage[all]{hypcap}       
 
\newcommand{\C}{\mathbb{C}}
\newcommand{\R}{\mathbb{R}}
\newcommand{\Q}{\mathbb{Q}}
\newcommand{\T}{\mathbb{T}}
\newcommand{\Z}{\mathbb{Z}}
\newcommand{\HH}{\mathbb{H}}
\newcommand{\opH}{\operatorname{H}}
\newcommand{\ii}{\operatorname{i}}
\newcommand{\eps}{\varepsilon}

\newcommand{\conv}{\operatorname{conv}}
\newcommand{\del}{\partial}
\newcommand{\delbar}{\bar{\partial}}
\newcommand{\PP}{\mathbb{P}}

\newcommand{\cA}{\mathcal{A}}
\newcommand{\cB}{\mathcal{B}}

\newcommand{\cX}{\mathcal{X}}
\newcommand{\cY}{\mathcal{Y}}

\newcommand{\cL}{\mathcal{L}}
\newcommand{\cM}{\mathcal{M}}
\newcommand{\cN}{\mathcal{N}}
\newcommand{\cP}{\mathcal{P}}
\newcommand{\cQ}{\mathcal{Q}}
\newcommand{\olo}{\mathcal{O}}
\newcommand{\cD}{\mathcal{D}}
\newcommand{\cT}{\mathcal{T}}
\newcommand{\cV}{\mathcal{V}}
\newcommand{\cZ}{\mathcal{Z}}

\newcommand{\barcT}{\overline{\cT}}
\newcommand{\Hom}{\operatorname{Hom}}
\newcommand{\Bl}{\operatorname{Bl}}

\newcommand{\res}{\operatorname{Res}}
\newcommand{\GM}{\operatorname{GM}}
\newcommand{\Jac}{\operatorname{Jac}}
\newcommand{\Crit}{\operatorname{Crit}}
\newcommand{\cW}{\mathcal{W}}
\newcommand{\Spec}{\operatorname{Spec}}
\newcommand{\Spf}{\operatorname{Spf}}
\newcommand{\kk}{\operatorname{k}}
\newcommand{\frk}{\mathfrak{k}}
\newcommand{\DF}{\operatorname{DF}}
\newcommand{\asym}{\operatorname{Asym}}
\newcommand{\bra}{\langle}
\newcommand{\ket}{\rangle}
\newcommand{\reg}{\operatorname{reg}}

\newcommand{\Bs}{\operatorname{Bs}}
\newcommand{\opT}{\operatorname{T}}
 
\newtheorem{thm}{Theorem}[section]
 
\newtheorem{prop}[thm]{Proposition}
\newtheorem{lemma}[thm]{Lemma}
\newtheorem{cor}[thm]{Corollary}

\theoremstyle{definition}
\newtheorem{definition}[thm]{Definition}

\theoremstyle{remark}
\newtheorem{exm}[thm]{Example}
\newtheorem{rmk}[thm]{Remark}

\title{Toric mirrors and test configurations}
\author{Jacopo Stoppa}
 
\begin{document}

\maketitle

\begin{abstract} We obtain results that relate Donaldson-Futaki type invariants (that is, the numerical invariants used to define $K$-stability for general polarised manifolds) for a toric polarised manifold and for a compactification of its mirror Landau-Ginzburg model, nearby the large volume limit. In general, these have the form of expansions containing terms which involve the base loci of certain linear systems determined by the Landau-Ginzburg potential (as expected from known constructions of compactified mirrors), and we give a condition under which these terms are subleading.

As an application we show that recently proposed notions of $K$-stability involving elements of the extended K\"ahler moduli space, i.e. $Z$-stability for polarised varieties, appear naturally from considerations of mirror symmetry (as a mirror to classical $K$-stability).
 \end{abstract}
\section{Introduction and main results}\label{Intro}
\subsection{Motivation from $K$-stability and mirror symmetry}
$K$-stability and mirror symmetry are both closely related to moduli theory: the former gives, conjecturally, the correct stability condition for constructing Hausdorff moduli spaces of polarised complex structures (thus involving at the same time complex and K\"ahler moduli), while the latter predicts that the complex and K\"ahler moduli of a mirror pair are intertwined. It seems natural to study the possible connections between the two theories.

The present paper is motivated by the problem of studying the behaviour of $K$-(semi)\-stability under mirror maps, in situations where (semi)stability is not automatic (while compact Calabi-Yau manifolds are automatically $K$-stable, see \cite{Odaka_CalabiKstab}). The most important example is given by Fano manifolds, which can be $K$-unstable, and whose mirror is expected to be a suitable Landau-Ginzburg model (see for example \cite{delcroix_report} for a recent collection of contributions aimed at emphasising the links between $K$-stability and Fano mirror symmetry). Here we work in the general \emph{toric} case, not necessarily Fano, where very strong mirror symmetry results are available. Note that even in the Fano case it is important to allow general polarisations, since these should correspond to the complex moduli of the mirror.

Our main results, Theorems \ref{MainThmIntro}, \ref{VanishingThm} and \ref{MainThmVariant}, give a precise relation between certain Donaldson-Futaki type invariants (the numerical invariants used to define $K$-stability) on the two sides of the mirror correspondence, at least nearby the large volume limit. To the author's knowledge, these are the first results relating stability invariants for mirror pairs (see Remark \ref{DualMomentumRmk}). Some conjectural expectations were proposed in \cite{ScarpaStoppa_complexified}, Section 1.7.1 and indeed our results provide partial confirmation for these (in particular, the appearance of twisted stability on the Landau-Ginzburg side).

Theorems \ref{MainThmIntro}, \ref{VanishingThm} show an important present limitation, since on the Landau-Ginzburg side, in general, we must replace the K\"ahler classes appearing in usual $K$-stability with \emph{complex} $(1,1)$-classes which, moreover, do not a priori satisfy \emph{semipositivity} conditions, and so the Donaldson-Futaki invariants on that side, although well-defined, do not have an immediate link to the usual notions of $K$-stability applied to the mirror. In our results,  reality and semipositivity must be imposed as further assumptions on the K\"ahler moduli (see Corollary \ref{PositiveCor}). 
\begin{rmk} In fact we will see that an even more natural result, Theorem \ref{MainThmVariant}, can be obtained if we allow the ``K\"ahler parameters" on the mirror to lie in the full \emph{extended K\"ahler moduli space} appearing in global versions of mirror symmetry (see our discussion in Section \ref{IntroAnaloguesSec}). This fits well with a recent proposal of generalised $K$-stability notions, due to Dervan, known as \emph{$Z$-stability for polarised varieties}, see Remark \ref{DervanRmk}. We believe that Theorem \ref{MainThmVariant} is conceptually relevant in this context, as it shows why $Z$-stability for polarised varieties appears naturally from considerations of mirror symmetry (as a mirror to classical $K$-stability).  
\end{rmk}
Let $(X, L)$ be a polarised toric manifold of dimension $n$. As we recall in Section \ref{BackSec}, it admits a Hodge-theoretic mirror given by a toric Landau-Ginzburg (LG) model, that is, roughly speaking, a nonconstant morphism 
\begin{equation*}
W_X(L)\!: \opT \to \C,\, \opT \cong (\C^*)^n,
\end{equation*}
known as the LG potential, where $\opT$ is dual to structure torus of $X$. 

As in the previous paper \cite{Stoppa_LargeComplex}, we are interested in understanding what conditions must be satisfied by the LG potential $W_X(L)$ if it is mirror to a $K$-semistable pair $(X, L)$ (the present paper can be read independently of \cite{Stoppa_LargeComplex}).

Recall that $K$-semistability is defined by requiring that a numerical weight attached to certain polarised one-parameter degenerations, known as the Donaldson-Futaki invariant of test configurations, is nonnegative (see \cite{Gabor_book} for a general introduction). 

If $(\cX, \cL) \to \PP^1$ is a regular compactified toric test configuration for $(X, L)$, i.e. the total space $\cX$ is smooth, toric and the map is equivariant, with $\cL$ ample (all this can be assumed for testing $K$-semistability), we can apply Hodge-theoretic mirror symmetry to the family $(\cX, \cL)$ (in the very general form appearing in the work \cite{CoatesCortiIritani_hodge} of Coates, Corti, Iritani and Tseng, which also contains an extensive list of references to previous results due to several authors). As we recall in Section \ref{DFResSec}, this allows to write the Donaldson-Futaki invariant $\DF(\cX, \cL)$ intrinsically on the mirror LG model 
\begin{equation*}
\cW_{\cX}(\cL)\!: \T \to \C,\,\T \cong (\C^*)^{n+1}, 
\end{equation*}
where $\T$ is dual to the structure torus of $\cX$. 

This involves the Grothendieck residue pairing of two specific (classes of) meromorphic forms $\alpha_{\cX}(\cL), \beta_{\cX}(\cL)$ in the top cohomology of the twisted de Rham complex $(\Omega^{\bullet}(\T), d\cW_{\cX}(\cL)\wedge)$, given explicitly in \eqref{alphaForm}, \eqref{betaForm}, such that 
\begin{equation*}
\DF(\cX, \cL) = \res_{\cW_{\cX}(\cL)}(\alpha_{\cX}(\cL), \beta_{\cX}(\cL)).
\end{equation*}
In the case of the trivial test configuration $X \times \PP^1$, with trivial action, the Donaldson-Futaki invariant vanishes of course, and abusing notation slightly we write this condition on the mirror as 
\begin{equation*}
\res_{W_{X}(L)}(\alpha_{X}(L), \beta_{X}(L)) = 0.
\end{equation*}
Moreover, as we explain in Section \ref{DFResSec}, it is always possible to regard $\cW_{\cX}(\cL)$ as a deformation of the original LG potential $W_X(L)$, in a precise sense, and similarly $\alpha_{\cX}(\cL), \beta_{\cX}(\cL)$ are deformations of $\alpha_{X}(L), \beta_{X}(L)$, so that \emph{$K$-semistability becomes the condition that, for all such deformations, the Grothendieck residue remains semipositive}, 
\begin{equation*}
\res_{\cW_{\cX}(\cL)}(\alpha_{\cX}(\cL), \beta_{\cX}(\cL)) \geq 0.
\end{equation*} 
\begin{exm} A basic, well-known $K$-destabilising example is given by taking $\cX$ to be the degeneration to the normal cone of the exceptional divisor $E \subset (X, L) := (\Bl_p \PP^2, -K_X)$ (see \cite{RossThomas_obstruction}, Example 5.27), $\cX := \Bl_{E \times \{0\}} X \times \PP^1$, $\cL : = -K_{\cX}$. In this case we have, in suitable torus coordinates,
\begin{align*}
& W_X(-K_X) = x + x' + x x' + \frac{1}{x x'},\\
&  \cW_{X \times \PP^1}(-K_{X \times \PP^1}) = \frac{1}{z} + \left(x + x' + x x' + \frac{1}{x x'}\right) + z,\\
& \cW_{\cX}(-K_{\cX}) = \frac{1}{z} + \left(x + x' + x x' + \frac{1}{x x'}\right) + z (1+ x x'),\\ 
& \res_{\cW_{\cX}(-K_{\cX})}(\alpha_{\cX}(\cL), \beta_{\cX}(\cL)) = \DF(\cX, \cL) < 0, 
\end{align*}
(see Example \ref{SlopeUnstableSurfaceExmPotentials})), showing a deformation of $W_X(-K_X)$ with negative Grothendieck residue.
\end{exm}
We would like to understand whether this semi-positivity condition can be translated into an algebro-geometric notion of \emph{stability} for the LG model $W_X(L)\!: (\C^*)^{n} \to \C$, or perhaps rather of a suitable \emph{compactification}. 
\subsection{Main results}
Our results in the present paper are stated in terms of a notion of formal Donaldson-Futaki invariant. Recall that with our assumptions the quantity $\DF(\cX, \cL)$ can be expressed through Atiyah-Bott localisation on $\cX$, with respect to the $\C^*$-action underlying the test configuration, as
\begin{align}\label{DFLocIntro}
\DF(\cX, \cL) = \sum_Z \int_Z\frac{(c_1(\cL) - h_{\cL})^n \left( \frac{n c}{n+1}(c_1(\cL) - h_{\cL}) + c_1(K_{\cX}) - \sum^{n-n_Z}_{i=0} w_i(v) + 1\right)}{e(N^{\cX}_Z(v))},
\end{align}
where $Z$ are the $\C^*$-fixed submanifolds, $c = \frac{c_1(X) \cup (c_1(L))^{n-1}}{(c_1(L))^n}$, $h_{\cL}$ is the Hamiltonian for the $S^1$-action, with generator denoted by $v$, with respect to (a Chern representative of) $c_1(\cL)$, $K_{\cX}$ is the canonical bundle, $w_i$ denote the weights of the $\C^*$-action on the normal bundle $N^{\cX}_Z(v)$ and $e(N^{\cX}_Z(v))$ is the equivariant Euler class (see in particular \cite{Legendre_localised}, Section 5). 
\begin{rmk}\label{TopConstRmk} By the definition of a compactified test configuration, a fixed submanifold $Z$ is either contained in the central fibre $\cX_0$ or equal to the fibre at infinity $\cX_{\infty}$. According to \cite{Legendre_localised}, Proposition 5.1, the choice of the topological constant $c$ is precisely such that the contribution from $\cX_{\infty}$ vanishes, i.e. the sum in \eqref{DFLocIntro} is effectively over $Z \subset \cX_0$. 
\end{rmk}
\begin{definition} If $\eta$ is a closed, \emph{complex} $(1,1)$-form on $\cX$, invariant under the $S^1$-action and with $\int_X \eta^n \neq 0$ but \emph{without positivity assumptions} (on its real and imaginary parts), and $h_{\eta}$ satisfies $\iota_v \eta = \delbar h_{\eta}$ for the generator $v$ of the $S^1$-action, then we \emph{define} the formal invariant $\DF(\cX, [\eta])$ as the right hand side of \eqref{DFLocIntro}, where $c_1(\cL)$, $h_{\cL}$ are replaced by $\eta$, $h_{\eta}$ (and so $c_1(L)$ is replaced by $[i^*_X \eta]$, where $i\!: X \cong \cX_1 \subset \cX$ embeds as the fibre $\cX_1$).
\end{definition}

Similarly, there is a notion of \emph{twisted} Donaldson-Futaki invariant $\DF(\cX, \cL, \cN)$ for test configurations of a triple $(X, L, N)$, where $(X, L)$ is a polarised manifold and $N$ is an additional holomorphic line bundle on $X$, introduced in \cite{Dervan_twistedcscK} (when $N$ is positive, this is the relevant notion for twisted constant scalar curvature K\"ahler metrics, with twist lying in $c_1(N)$, see e.g. \cite{Dervan_twistedcscK, StoppaTwisted}). This satisfies the localisation formula 
\begin{align}\label{DFMapsLocIntro}
&\nonumber\DF(\cX, \cL,\cN)\\ 
&= \sum_Z \int_Z\frac{(c_1(\cL) - h_{\cL})^n \left( \frac{n c_N}{n+1}(c_1(\cL) - h_{\cL}) + c_1(K_{\cX}) - \sum^{n-n_Z}_{i=0} w_i(v) + 1 + c_1(\cN) - h_{\cN}\right)}{e(N^{\cX}_Z(v))},
\end{align}
where $h_{\cN}$ is Hamiltonian with respect to $c_1(\cN)$ and $c_N = \frac{(c_1(X) - c_1(N))\cup (c_1(L))^{n-1}}{(c_1(L))^n}$. 

\begin{definition} We extend \eqref{DFMapsLocIntro} formally to $\DF(\cX, [\eta], [\xi])$, using the right hand side as a \emph{definition}, where $c_1(\cL)$, $c_1(\cN)$ are replaced by complex $(1, 1)$-forms $\eta$, $\xi$, with complex Hamiltonians $h_{\eta}$,$h_{\xi}$. 
\end{definition}
Given this, we can state our first main result, which concerns the \emph{large volume limit}, when $L$ is replaced by $k L := L^{\otimes k}$ for $k \gg 1$. It is conditional on a property of a certain cohomology intersection pairing attached to the LG potential, known as the stationary phase formula, which we spell out in Section \ref{DFchSubSec}.  
\begin{thm}\label{MainThmIntro} Fix a regular toric test configuration $(\cX, \cL)$ for the polarised toric manifold $(X, L)$. There exists a compactification $\barcT$ of the torus $\T$, endowed with a logarithmic connection $\nabla_{\cW_{\cX}(k\cL)}$ determined by the mirror Landau-Ginzburg potential $\cW_{\cX}(k\cL)$, with the property that, if the stationary phase formula holds for $(\barcT, \nabla_{\cW_{\cX}(k\cL)})$ for $k\gg 1$ (see condition $\dagger$ in Corollary \ref{FutakiToCohoIFCor}), then one can construct regular twisted toric test configurations for a compactification $X^{\vee}$ of the torus $\opT$,  
\begin{equation*}
(\cX^{\vee, i}, [\eta^i_k], [\xi^i_k]),\,i = 1, \ldots, m, 
\end{equation*}
for complex $(1,1)$-classes $[\eta^i_k]$, $[\xi^i_k]$ (without semipositivity conditions), such that their formal Donaldson-Futaki invariants are well-defined and satisfy
\begin{align*}
&\DF(\cX, \cL) = \res_{\cW_{\cX}(k\cL)}(k^{-n}\alpha_{\cX}(k\cL), \beta_{\cX}(k\cL))\\
& = \sum^m_{i=1}\DF(\cX^{\vee, i}, [\eta^i_k], [\xi^i_k]) + \cB(\cW_{\cX}(k\cL))+ O(k^{-1}),
\end{align*}
where $\cB(\cW_{\cX}(k\cL))$ is an explicit quantity determined by the (non-toric) base loci of certain pencils defined by $\cW_{\cX}(k\cL)$ on an intermediate toric compactification $\barcT_{\T}$, with $\T \subset \barcT_{\T} \subset \barcT$ (see \eqref{BaseLocusTerm}).
\end{thm}
The proof of Theorem \ref{MainThmIntro} is completed in Section \ref{MainProofSubSecCompletion}, after setting up the construction in  Sections \ref{DFchSubSec}, \ref{MainProofSubSec}, using certain facts from the theory of logarithmic connections (recalled in Section \ref{LogBackSec}). It relies heavily on some aspects of toric Hodge-theoretic mirror symmetry after Coates, Corti, Iritani and Tseng \cite{CoatesCortiIritani_hodge}, and these are briefly recalled in Section \ref{BackSec}.
\begin{rmk}\label{MainThmRmk} Theorem \ref{MainThmIntro} requires some comments.
\begin{itemize}
\item[$(i)$] The compactification $X^{\vee}$, the total spaces $\cX^{\vee, i}$, the ``polarisations" $[\eta^i_k]$ and the twists $[\xi^i_k]$ are determined by the LG potentials $W_X(L)$, $\cW_{\cX}(\cL)$ or, more precisely, by the pair $(\barcT, \nabla_{\cW_{\cX}(k\cL)})$. However, $X^{\vee}$, $m$ and the total spaces $\cX^{\vee, i}$, $i = 1, \ldots, m$ can be chosen independently of $k$. The construction of the total spaces $\cX^{\vee, i}$ from $\barcT$ is described in Section \ref{MainProofSubSec}, and that of the $(1, 1)$-classes $[\eta^i_k]$, $[\xi^i_k]$ in Section \ref{MainProofSubSecCompletion}.

In particular, we will see that when the general fibre $X := T_{P}$ is a toric Fano with reflexive fan polytope $P$, then $X^{\vee}$ can be chosen as a resolution of the toric variety $T_{P^{\circ}}$ given by the polar dual $P^{\circ}$, as usual in Fano toric mirror symmetry.
\item[$(ii)$] The appearance of the term $\cB(\cW_{\cX}(k\cL))$ (described in \eqref{BaseLocusTerm}) is expected: it is known that, at least in the case when $\cX = T_{\cP}$ is a Fano surface or threefold with reflexive polytope $\cP$, a genuine compactified mirror can be obtained from the polar dual $T_{\cP^{\circ}}$ by resolving singularities and blowing up the base locus of the pencil defined by $\cW_{\cX}(\cL)$ and the toric boundary (see e.g. \cite{Przyjalkowski_CYLG}). In this case, the contributions $\DF(\cX^{\vee, i}, [\eta^i_k], [\xi^i_k])$ are computed on a \emph{toric} resolution of $T_{\cP^{\circ}}$, while $\cB(\cW_{\cX}(k\cL))$ represents the contribution from the \emph{non-toric} blowups needed to resolve the relevant base loci.
 
\item[$(iii)$] As will be clear from the proof, since we are using the mirror map, the subleading term $O(k^{-1})$ is determined by the quantum cohomology of $\cX$ (which degenerates to usual intersection theory as $k \to \infty$), in a complicated way. It would be interesting to see if the subleading term can be understood as a suitable deformation of the invariants $\DF(\cX^{\vee, i}, [\eta^i_k], [\xi^i_k])$.
\item[$(iv)$] It is not hard to recognise that the quantity $\res_{\cW_{\cX}(k\cL)}(k^{-n}\alpha_{\cX}(k\cL), \beta_{\cX}(k\cL))$, for $k \gg 1$, resembles a localised Donaldson-Futaki invariant, where, however, the total space is the affine torus $\T$, and the holomorphic vector field is given by $\nabla \cW(k\cL)$ (see Remark \ref{NoncompactDFRmk}). The idea of Theorem \ref{MainThmIntro} is to use a suitable residue theorem in order to replace $(\C^*)^{n+1}$ with compactified test configurations. This requires relating $\res_{\cW_{\cX}(k\cL)}(k^{-n}\alpha_{\cX}(k\cL), \beta_{\cX}(k\cL))$, which is defined by a ``stationary phase formula" (see Section \ref{HigherResSec}), to a global intersection number. This is the role of the stationary phase formula appearing in Theorem \ref{MainThmIntro}. 
\item[$(v)$] Understanding when the stationary phase formula holds is an important problem studied in different contexts, see e.g. the discussion in \cite{Matsubara_cohomology}, Section 2 (briefly recalled in our Section \ref{DFchSubSec}), which relies on recent results of Huh \cite{huh_likelihood}.
\item[$(vi)$] The classic work of Donaldson \cite{Donaldson_stability_toric} gives a complete description of toric test configurations and their Donaldson-Futaki invariants using convex geometry on the momentum polytope (see \cite{Yanir_higher} for a recent extension of these results to other stability invariants). However it seems that this description is not well suited for the application of toric mirror theorems, see Section \ref{BackSec}.
\item[$(vii)$] The famous conjecture that $K$-polystability is equivalent to the existence of constant scalar curvature K\"ahler metrics is known for toric surfaces \cite{Donaldson_cscKmetrics_toricsurfaces} and the toric uniform case is known in all dimensions \cite{ChenCheng_existence}.  
\end{itemize}
\end{rmk}
In order to relate $K$-(semi)stability more directly on the two sides of the mirror correspondence, it is natural to ask when the term $\cB_k(\cW_{\cX}(\cL))$, determined by the base loci, is actually subleading, i.e. we have $\cB_k(\cW_{\cX}(\cL)) = O(k^{-1})$. In the light of Remark \ref{MainThmRmk} $(ii)$, we might call this the Arezzo-Pacard limit, in which the Donaldson-Futaki invariant of a test configuration for a blowup is dominated by that of the base (see e.g. \cite{Stoppa_JAG}). From this viewpoint, at least heuristically, we are asking what K\"ahler classes on $\cX$, nearby the large volume limit, correspond to the Arezzo-Pacard limit on the mirror, i.e. for $(\cX^{\vee, i}, [\eta^i_k], [\xi^i_k])$ and $k \gg 1$. 
\begin{thm}\label{VanishingThm} In the setup of Theorem \ref{MainThmIntro}, let $\cA$ be an ample line bundle on $\cX$. Suppose that, in terms of the linear mirror isomorphism $\Theta$ at $k \cA$ (see Section \ref{BackSec}), for some $r > 0$, we have 
\begin{equation*}
\cW_{\cX}(k \cA) = r^{-1}\Theta_{k \cA }(\cA) + O(k^{-1}).  
\end{equation*}
This is the case for example if $\cX$ is Fano and $\cA = - r K_{\cX}$ is a positive multiple of $-K_{\cX}$, so by properties of $\Theta$ we have $\cW_{\cX}(k(- r K_{\cX})) = \Theta_{k(- r K_{\cX})}(-K_{\cX}) = r^{-1} \Theta_{k(- r K_{\cX})}(- r K_{\cX})$ (see \eqref{FanoMirrorOfAnticanonical}). 
 
Then, we have $\cB_k(\cW_{\cX}(k \cA )) = O(k^{-1})$, so
\begin{align*}
\DF(\cX, \cA) = \sum^m_{i=1}\DF( \cX^{\vee, i}, [\eta^i_k], [\xi^i_k]) + O(k^{-1}).
\end{align*}
\end{thm} 
Theorem \ref{VanishingThm} is proved in Section \ref{MainProofSubSecCompletion}.
\begin{exm} Let $(\cX := \Bl_{p \times \{0\}} \PP^1 \times \PP^1, -\frac{1}{2}K_{\cX})$ be the test configuration for $(X, L) = (\PP^1, \olo(1))$ given by the degeneration to the normal cone of a point. We will show (see Examples \ref{NormalConeAnticanTest}, \ref{NormalConeAnticanExm}) that in this case the stationary phase formula holds, one can take $X^{\vee} \cong \PP^1$, $\cX^{\vee, 1} \cong \Bl_{q_1 \times \{0\}, q_2 \times \{0\}} \PP^1 \times \PP^1$ for $q_1 \neq q_2$, $\cX^{\vee, 2} \cong \Bl_{q\times\{0\}} \PP^1 \times \PP^1$ endowed with suitable K\"ahler classes $[\eta^1_{k}]$, $[\eta^2_k]$ and twists $[\xi^1_k]$, $[\xi^2_k]$, and we have
\begin{align*}
& \frac{1}{4} = \DF(\cX, -\frac{1}{2} K_{\cX}) = \DF(\cX^{\vee, 1}, [\eta^1_k], [\xi^1_k]) + \DF(\cX^{\vee, 2}, [\eta^2_k], [\xi^2_k]) + O(k^{-1}),\\  
& \DF(\cX^{\vee, 1}, [\eta^1_k], [\xi^1_k]) = O(k^{-1}),\,\DF(\cX^{\vee, 2}, [\eta^2_k], [\xi^2_k]) = \frac{1}{4} + O(k^{-1}).
\end{align*}
It is known that, for smooth curves, it is sufficient to test $K$-polystability on degeneration to the normal cone of a single point (see \cite{RossThomas_obstruction}, Section 5.3). So in our example it is sufficient to test on $(\cX, -\frac{1}{2}K_{\cX})$, and this is self-mirror in the sense that $\cX \cong \cX^{\vee, 2}$ and $\DF(\cX, -\frac{1}{2} K_{\cX}) = \DF(\cX^{\vee, 2}, [\eta^2_k], [\xi^2_k])$ up to $O(k^{-1})$.
\end{exm}
\begin{exm} Consider $(\cX := \PP(\olo \oplus \olo(1)), -\frac{1}{3}K_{\cX})$ regarded as a product test configuration for $\PP^1$. It is also known that for $\PP^1$ it is sufficient to test $K$-polystability on this product test configuration (e.g. by a very special case of the results of \cite{WangZhu_toric}). We will show (see Examples \ref{ProductAnticanTest}, \ref{ProductAnticanExm}) that in this case the stationary phase formula holds, one can take $X^{\vee} \cong \PP^1$, $\cX^{\vee, i} \cong \Bl_{q' \times \{0\}, q''\times\{0\}} \PP^1 \times \PP^1$ for $i =1, 2$ (with $q' \times \{0\}, q''\times\{0\}$ infinitely close), endowed with suitable K\"ahler classes $[\eta^i_{k}]$ and twists $[\xi^i_k]$, and we have
\begin{align*}
& 0 = \DF(\cX, -\frac{1}{2} K_{\cX}) = \DF( \cX^{\vee, 1}, [\eta^1_k], [\xi^1_k]) + \DF(\cX^{\vee, 2}, [\eta^2_k], [\xi^2_k]) + O(k^{-1}),\\  
& \DF(\cX^{\vee, i}, [\eta^i_k], [\xi^i_k]) = O(k^{-1}),\, i =1, 2.
\end{align*}
\end{exm}
Theorem \ref{VanishingThm} and Remark \ref{MainThmRmk} $(i)$ have the following obvious consequences.  
\begin{cor}\label{PositiveCor} In the situation of Theorem \ref{VanishingThm}, suppose the K\"ahler parameter $\cA$ on $\cX$ is such that all the classes $[\eta^i_k]$, $[\xi^i_k]$ are real and semipositive for $k \gg 1$ and we have $\DF(\cX, \cA) < 0$ (so $(X, A:= \cA|_{X})$ is $K$-unstable). Then, $(X^{\vee}, [\eta^i_k]|_{X^{\vee}}, [\xi^i_k]|_{X^{\vee}})$ is twisted $K$-unstable for some $i$ and sufficiently large $k$. In general, with the assumptions of Theorem \ref{VanishingThm}, $\DF(\cX, \cA) < 0$ implies that $ \Re \DF( \cX^{\vee, i}, [\eta^i_k], [\xi^i_k]) < 0$ for some $i$ and sufficiently large $k$, but $[\eta^i_k]$, $[\xi^i_k]$ are not necessarily real and semipositive. 
\end{cor}
\begin{cor} In the situation of Theorem \ref{VanishingThm}, suppose that $(X, A:=\cA|_X)$ is $K$-semistable. Then, for all $(\cX, \cA)$ with corresponding $(\barcT, \nabla_{\cW_{\cX}(k\cA)})$, we have 
\begin{equation*}
\liminf_{k \to \infty} \max_{i} \Re \DF( \cX^{\vee, i}, [\eta^i_k], [\xi^i_k]) \geq 0. 
\end{equation*}
\end{cor}
Naturally, analogous statements hold in the general case of Theorem \ref{MainThmIntro}, but they involve the term $\cB_k(\cW_{\cX}(\cL))$, that is, we have
\begin{equation*}
\min \{\{\Re \DF(\cX^{\vee, i}, [\eta^i_k], [\xi^i_k])\}_i, \Re\cB_k(\cW_{\cX}(\cL))\} < 0 \textrm{ for } k\gg 1,
\end{equation*}
respectively
\begin{equation*}
\liminf_{k \to \infty} \max \{\{\Re \DF( \cX^{\vee, i}, [\eta^i_k], [\xi^i_k])\}_i, \Re\cB_k(\cW_{\cX}(\cL))\} \geq 0. 
\end{equation*}
\begin{exm} Corollary \ref{PositiveCor} can be applied to the example of degeneration to the normal cone $\cX := \Bl_{E \times \{0\}} X \times \PP^1$, $\cL : = -K_{\cX}$ of the exceptional divisor $E$ which $K$-destabilises $(X, L) = (\Bl_p \PP^2, -K_{\Bl_p}\PP^2)$. In this case, $X^{\vee}$ is the crepant resolution of the orbifold $T_{P^{\circ}}$ given by a toric, singular intersection of two quadrics in $\PP^4$. Assuming the stationary phase formula, Theorem \ref{VanishingThm} provides a (non-positive)``destabilising" twisted test configuration $(\cX^{\vee}, [\eta], [\xi])$ for $X^{\vee}$, \emph{which is mirror to degeneration to the normal cone of $E$}, and this can be made more explicit if we work with the orbifold $T_{P^{\circ}}$ itself, see Examples \ref{SlopeUnstableSurfaceExm}, \ref{SlopeUnstableSurfaceExmVanishing}. Note that it is known that both $\Bl_p \PP^2$ and the orbifold $T_{P^{\circ}}$ are $K$-unstable with respect to all K\"ahler classes, since their automorphisms groups are non-reductive. In our construction, however, there is a complex twist $[\xi]$ and the ``K\"ahler class" $[\eta]$ is a priori only a complex $(1, 1)$-class; both are determined by the LG potential $\cW_{\cX}(\cL)$.
\end{exm}

\begin{rmk}\label{DualMomentumRmk} An anticanonically polarised Fano  $(X, -K_X)$ is $K$-polystable iff the classical Futaki character on holomorphic vector fields, identified with the barycenter $b(P^{\circ})$ of the momentum polytope $P^{\circ}$ (the polar dual of the fan polytope $P$) of $(X, -K_X)$, vanishes (see \cite{WangZhu_toric}). Thus, in view of polar duality as a first approximation to Fano toric mirror symmetry (see Remark \ref{MainThmRmk}, $(ii)$), it is natural to study the relation between the barycenters $b(P^{\circ})$, $b(P)$. Simple examples show that the condition $b(P) = 0$ does \emph{not} imply $b(P^{\circ}) = 0$, except in the very special symmetric case when $P = - P$.   

Sano \cite{Sano_polar} (with motivations not coming from mirror symmetry) defines and studies a ``polar dual" $e_P \in P$ for a suitable multiple of the barycenter of $P^{\circ}$, and shows that in low dimensions the vanishing of $e_P$ implies the condition $b(P^{\circ}) = 0$. However this implication breaks down in higher dimensions. From our perspective, the correct approach to this problem is to regard it as a special case of that of relating the Donaldson-Futaki invariants of $X$, $X^{\vee}$, when we start from a product test configuration for $(X, -K_X)$ (i.e. one induced from a holomorphic vector field); our results Theorems \ref{MainThmIntro}, \ref{MainThmVariant} then yield a relation between the barycenter $b(P^{\circ})$ and certain Donaldson-Futaki invariants computed on $X^{\vee}$. 
\end{rmk}

\subsection{Analogues involving the extended K\"ahler moduli space}\label{IntroAnaloguesSec} 
As will be clear from the proof of Theorem \ref{MainThmIntro}, the need to allow a potentially large set of test configurations $(\cX^{\vee, i}, [\eta^i_k], [\xi^i_k]),\,i = 1, \ldots, m$ on the compactified LG side arises because \emph{we insisted that the classes $[\eta^i_k]$, $[\xi^i_k]$ should be of degree $(1,1)$}, thus retaining some geometric interpretation as potential (complefixied) K\"ahler classes (or complex twist classes) on the mirror. 

It turns out that even more natural variants of Theorems \ref{MainThmIntro}, \ref{VanishingThm} can be obtained, essentially with the same proofs, if we allow $[\eta], [\xi] \in H^{1, 1}(\cX^{\vee}, \C)$ to be replaced with cohomology classes $\cZ_k, \tilde{\cZ}_k \in H^*(\cX^{\vee}, \C)$ \emph{of arbitrary degrees}. 

While this replacement seems far from classical notions in $K$-stability (see Remark \ref{DervanRmk}), it appears natural from the mirror symmetry perspective, where $H^*(\cX^{\vee}, \C)$ plays the role of the \emph{extended K\"ahler moduli space of $\cX^{\vee}$} (this is a well-known general feature of \emph{global} mirror symmetry, and in the toric case it is explained and studied in depth in \cite{IritaniGlobalMirror}).

\begin{definition}\label{DFExtendedKahler} Suppose $(\cX^{\vee}, v)$ is the total space of a regular, compactified toric test configuration for the fibre $X^{\vee} \cong \cX^{\vee}_1$, with $\C^*$-action generator $v$. Fix any pair of classes
$\cZ, \tilde{\cZ} \in H^*(\cX^{\vee}, \C)$ (``extended K\"ahler parameters").
We define the Donaldson-Futaki invariant of $(\cX^{\vee}, v, \cZ , \tilde{\cZ})$, regarded as a test configuration for the general fibre, i.e. for $(X^{\vee} := \cX^{\vee}_1, \cZ|_{X^{\vee}}, \tilde{\cZ}|_{X^{\vee}})$, as the right hand side of \eqref{DFMapsLocIntro} where the equivariant representatives of $[\eta], [\xi]$ are replaced with equivariant representatives $\cZ' , \tilde{\cZ}' \in H^*_{S^1}(\cX^{\vee}, \C)$ of $\cZ$, $\tilde{\cZ}$, namely  
\begin{align*} 
\DF(\cX^{\vee}, \cZ, \tilde{\cZ}) = \sum_Z \int_Z\frac{(\cZ')^n \left( \frac{n}{n+1} c_{\cZ, \tilde{\cZ}'} \cZ' + c_1(K_{\cX}) - \sum^{n-n_Z}_{i=0} w_i(v) + 1 + \tilde{\cZ}'\right)}{e(N^{\cX}_Z(v))},
\end{align*}
while the constant $c_{\cZ, \tilde{\cZ}'} \in \C$ is chosen so that the contribution from the fixed component $\cX^{\vee}_{\infty}$ vanishes (as in Remark \ref{TopConstRmk}).
\end{definition}
\begin{thm}\label{MainThmVariant} Suppose the same assumptions of Theorem \ref{MainThmIntro} hold. Then there exist regular twisted toric test configurations for a compactification $X^{\vee}$ of the torus $\opT$,  
\begin{equation*}
(\cX^{\vee, 0}, \cZ^0_k, \tilde{\cZ}^0_k),\,(\cX^{\vee, \infty}, \cZ^{\infty}_k, \tilde{\cZ}^{\infty}_k),  
\end{equation*}
$\cZ^0_k, \tilde{\cZ}^0_k \in H^*(\cX^{\vee, 0}, \C)$, $\cZ^{\infty}_k, \tilde{\cZ}^{\infty}_k \in H^*(\cX^{\vee, \infty}, \C)$ denoting extended K\"ahler parameters, with fixed restrictions to $X^{\vee}$, such that their Donaldson-Futaki invariants in the sense of Definition \ref{DFExtendedKahler} are well-defined and satisfy
\begin{align*}
\DF(\cX, \cL) = \DF(\cX^{\vee, 0}, \cZ^0_k, \tilde{\cZ}^0_k) + \DF(\cX^{\vee, \infty}, \cZ^{\infty}_k, \tilde{\cZ}^{\infty}_k) + \cB(\cW_{\cX}(k\cL))+ O(k^{-1}),
\end{align*}
where $\cB(\cW_{\cX}(k\cL))$ is the same contribution from non-toric base loci appearing in Theorem \ref{MainThmIntro}. Moreover, the splitting into $\cX^{\vee, 0}$, $\cX^{\vee, \infty}$ is given naturally by the stationary phase formula, see Definition \ref{TestConfigDef}. 

Similarly, under the additional assumptions of Theorem \ref{VanishingThm}, we have
\begin{align}\label{MainThmVariantVanishing}
\DF(\cX, \cA) = \DF(\cX^{\vee, 0}, \cZ^0_k, \tilde{\cZ}^0_k) + \DF(\cX^{\vee, \infty}, \cZ^{\infty}_k, \tilde{\cZ}^{\infty}_k) + O(k^{-1}).
\end{align} 
\end{thm}
Theorem \ref{MainThmVariant} is proved in Section \ref{MainProofSubSecCompletion}.

\begin{rmk}\label{DervanRmk} Dervan \cite{Dervan_crit_metrics} develops a notion of $K$-stability (known as \emph{$Z$-stability for polarised varieties}), in which, for any test configuration $(\cX, \cL)$ for a polarised variety $(X, L)$, an element of the extended K\"ahler moduli space $\cZ \in H^*(\cX, \C)$ is used to define the relevant Donaldson-Futaki invariant (in $\cite{Dervan_crit_metrics}$, $\cZ$ is regarded as a \emph{central charge}, i.e. dually as a map $H^*(\cX, \C) \to \C$, via pairing with the Chern character of the polarisation $\cL$). Although there are some technical differences (i.e. Dervan's choice of $\cZ$ must satisfy some constraints), our Theorem \ref{DFExtendedKahler} clearly shows why $Z$-stability for polarised varieties appears naturally from considerations of mirror symmetry, as a mirror to classical $K$-stability (see also \cite{Corradini} for closely related localisation formulae).       
\end{rmk}

\begin{exm} The case of $(\cX := \Bl_{p \times \{0\}} \PP^1 \times \PP^1, -\frac{1}{2}K_{\cX})$ discussed above also gives an example of \eqref{MainThmVariantVanishing} in Theorem \ref{MainThmVariant}, with $\cX^{\vee, 0} \cong \Bl_{q_1 \times \{0\}, q_2 \times \{0\}} \PP^1 \times \PP^1$ for $q_1 \neq q_2$, $\cX^{\vee, \infty} \cong \Bl_{q\times\{0\}} \PP^1 \times \PP^1$, $\cZ^0_k = [\eta^1_{k}]$, $\tilde{\cZ}^{0}_k = [\xi^1_k]$, $\cZ^{\infty}_k = [\eta^2_{k}]$, $\tilde{\cZ}^{\infty}_k = [\xi^2_k]$. 

Similarly, $(\cX := \PP(\olo \oplus \olo(1)), -\frac{1}{3}K_{\cX})$ gives an example of \eqref{MainThmVariantVanishing} with $\cX^{\vee, 0} \cong \cX^{\vee, \infty} \cong \Bl_{q' \times \{0\}, q''\times\{0\}} \PP^1 \times \PP^1$.
\end{exm}

\noindent{\textbf{Acknowledgements.}} I am grateful to Ruadha\'i Dervan, Hiroshi Iritani, Yuji Odaka, Yanir Rubinstein and Richard Thomas for helpful discussions related to the present work, and to the anonymous Referees for many important suggestions on improving the manuscript. This research was carried out in the framework of the project PRIN 2022BTA242 ``Geometry of algebraic structures: moduli, invariants, deformations''. 

\section{Some background}\label{BackSec}

\subsection{Toric LG models}\label{ToricLGSec} Let $X$ be a projective toric manifold. According to the general results of \cite{CoatesCortiIritani_hodge}, the mirror to $X$ is given by a Laurent polynomial of the form
\begin{equation*}
W(x; y) = \sum^m_{i = 1} y_i Q^{\lambda(b_i)} x^{b_i} + \sum_{ \kk \in G} y_{\kk} Q^{\lambda(\kk)} x^{\kk}, 
\end{equation*}
where $b_1, \ldots, b_m$ are generators of the rays of the fan of $X$ in its fan lattice $N$, the formal variable $Q$ is known as the Novikov variable, $\lambda(b_i), \lambda(\kk) \in H_2(X, \Q)$ are certain curve classes (see \cite{CoatesCortiIritani_hodge}, Section 4.1), $x \in \Hom(N, \C^*)$ is a torus coordinate, $y_i, y_{\kk}$ are deformation parameters, and $G \subset N$ is a suitable finite subset. Thus, fixing all parameters except $x$, $W(x; y)$ defines a regular function on the algebraic torus $\Hom(N, \C^*)$, known as the Landau-Ginzburg (LG) potential. 

More globally, $W$ must be regarded as a morphism $W\!: \widehat{\cY} \to \C$ where $\widehat{\cY} \to \widehat{\cM}$ is a degenerating family of affine toric varieties over the base $\widehat{\cM} = \Spf \C[\![\Lambda_{+}]\!] \times \Spf \C[\![y]\!]$ and $\Lambda_+ \subset H_2(X, \Q)$ denotes the monoid of effective curve classes, so $\C[\![\Lambda_{+}]\!]$ is the Novikov ring. There is a locus $\cM \subset \widehat{\cM}$ over which the family $\widehat{\cY} \to \widehat{\cM}$ restricts to a trivial fibration by algebraic tori $\cY \to \cM$, and a choice of trivialisation allows to write it as a Laurent polynomial depending on parameters, as above. A polarisation $L$ on $X$, or more generally a K\"ahler class $[\omega_X]$, defines a closed point of the base $\Spf \C[\![\Lambda_{+}]\!]$, and so specifies the exponents of the Novikov variable. 

Note that the extra term $\sum_{ \kk \in G} y_{\kk} Q^{\lambda(\kk)} x^{\kk}$ in the expression for $W(x; y)$ is necessary in the general case when $X$ is not a weak Fano manifold, and is required so that the parametrisation of $W(x; y)$ by K\"ahler classes of $X$ along a locus is effective. When $X$ is weak Fano, on the other hand, one can choose $G = \emptyset$. 

The Hodge-theoretic information attached to the LG mirror family is given by the \emph{logarithmic twisted de Rham complex}
\begin{equation*}
\left(\Omega^{\bullet}_{\widehat{\cY}/\widehat{\cM}}, z d+ d W\wedge\right).
\end{equation*}
The \emph{Gauss-Manin system} $\GM(W)$ is the top cohomology of this complex. It is naturally endowed with extra structure, in particular the \emph{Gauss-Manin connection} $\nabla^{\GM}$ (see \cite{CoatesCortiIritani_hodge}, Section 4.4) and the \emph{higher residue pairing} 
\begin{equation*}
P\!: \GM(W) \times \GM(W) \to \C[z] [\![\Lambda_+]\!][\![y]\!].
\end{equation*}
(discussed in \cite{CoatesCortiIritani_hodge}, Section 6). The latter is especially important for our purposes and will be recalled in Section \ref{HigherResSec}.

The main results of \cite{CoatesCortiIritani_hodge} (namely \cite{CoatesCortiIritani_hodge}, Theorems 4.28 and 6.11, summarised in \cite{CoatesCortiIritani_hodge}, Theorem 1.1) prove \emph{Hodge-theoretic mirror symmetry for the big quantum cohomology} in this case, part of which we recall briefly as follows:
\begin{itemize}
\item[$(i)$] the Gauss-Manin system is a free module over $\C[z] [\![\Lambda_+]\!][\![y]\!]$ of rank $\dim H^*(X, \C)$;
\item[$(ii)$] there are a mirror map $\tau = \tau(y) \in H^*(X, \C) \otimes \C [\![\Lambda_+]\!][\![y]\!]$ and a $\C[z] [\![\Lambda_+]\!][\![y]\!]$-linear mirror isomorphism
\begin{equation*}
\Theta\!: \GM(W) \xrightarrow{\,\,\cong\,\,} H^*(X, \C) \otimes \C[z] [\![\Lambda_+]\!][\![y]\!] 
\end{equation*}
such that  $\Theta$ intertwines the Gauss-Manin connection $\nabla^{\GM}$ with the pullback $\tau^*\nabla^{\operatorname{D}}$ of the quantum connection $\nabla^{\operatorname{D}}$ by $\tau$ (see \cite{CoatesCortiIritani_hodge} Section 3.2 for the latter). 

In the $z\to 0$ limit, choosing the generator $\Theta^{-1}(1)$ for the $\GM(W(L))|_{z=0}$ induces an algebra structure on $\GM( W(L))|_{z=0}$ and a linear mirror isomorphism 
\begin{equation*}
\Jac(W(L)) \cong (H^*(X, \C), *_{[c_1(L)]})
\end{equation*}
with the quantum cohomology ring, where we write $*_{[c_1(L)]}$ for the quantum product evaluated at the quantum parameter $Q$, such that
\begin{equation*} 
Q^d = \exp\left(-2 \pi \int_d c_1(L)\right),\,d \in \Lambda_{+}.
\end{equation*}
In the large volume limit $L \mapsto k L$ for $k \gg 1$, the linear mirror isomorphism satisfies
\begin{equation}\label{FanoMirrorOfAnticanonical}
\Theta_{kL}(W(k L)) = c_1(-K_{X}) + O(k^{-1}).
\end{equation}
\item[$(iii)$] The linear mirror isomorphism $\Theta$ intertwines the higher residue pairing on $\GM(W)$ with the Poincar\'e pairing on $H^*(X, \C) \otimes \C[z] [\![\Lambda_+]\!][\![y]\!]$, that is we have 
\begin{equation*}
P(\Omega_1, \Omega_2) = \left(\Theta(\Omega_1)\big|_{z \mapsto -z}, \Theta(\Omega_2)\right).
\end{equation*}
\end{itemize}
Note that $(i)$ includes a convergence property for the mirror map $\tau$ which we will recall in a moment. We can summarise $(i)$, $(ii)$, $(iii)$ roughly by saying that, for a fixed trivialisation of the torus fibration $\cY \to \cM$, a choice of K\"ahler class $[\omega_X]$ fixes the exponents of the Novikov variable in $W(x; y)$, but not the coefficients $y$. Then, as shown in \cite{CoatesCortiIritani_hodge}, Section 7, the formal power series mirror map $\tau$ is convergent at least after rescaling $[\omega_X]$ by a sufficiently large factor (this will be enough for our purposes), and it fixes the deformation parameters $y$ so that, for the corresponding LG potential $W(x; y) = W(\omega_X)$, the linear mirror isomorphism   
\begin{equation*}
\Theta_{\omega_X}\!: \GM(W(\omega_X)) \xrightarrow{\,\,\cong\,\,} H^*(X, \C) \otimes \C[\![z]\!] 
\end{equation*}
intertwines the Gauss-Manin and quantum connections, and preserves the higher residue pairing: 
\begin{equation*}
P_{W(\omega_X)}(\Omega_1, \Omega_2) = \left(\Theta_{\omega_X}(\Omega_1)\big|_{z \mapsto -z}, \Theta_{\omega_X}(\Omega_2)\right).
\end{equation*}
\subsection{Higher residue pairing}\label{HigherResSec}
It will be useful to know a little more concerning the higher residue pairing used in \cite{CoatesCortiIritani_hodge}. This is defined through a formal stationary phase expansion. Let us briefly review the construction. 

We start with a holomorphic function $f$ on $\C^n$ with a non-degenerate critical point $p$. Fix a stable manifold $\Gamma(p)$ for the Morse function $\Re(f(t))$ and perform the usual stationary phase expansion as $z \to 0$ for the integral
\begin{equation*}
\int_{\Gamma(p)} e^{f(t)/z} g(t) dt^1\cdots dt^n
\end{equation*} 
where $\Re(z) < 0$ and $g(t)$ is holomorphic. We write the result of the stationary phase expansion in the form
\begin{equation*}
e^{f(p)/z}(-2\pi z)^{n/2} \asym_p\left(e^{f(t)/z} g(t) dt\right). 
\end{equation*}
This is always possible and defines the formal power series $\asym_p\left(e^{f(t)/z} g(t) dt\right)$, which takes the form
\begin{equation*}
\frac{1}{\sqrt{\det(f_{ij}(p))}}\left(g(p) + a_1 z + a_2 z^2 + \cdots \right).
\end{equation*} 

Assume, for simplicity, that the LG potential $\cW_{\cX}(\cL)$ has nondegenerate isolated critical points. One shows that the formal power series $\asym_p\left(e^{\cW(\cL)/z} \Omega \right)$ vanishes if $e^{\cW(\cL)/z} \Omega$ is exact in the complex $\left(\Omega^{\bullet}, z d+ d \cW(\cL)\wedge\right)$, so the map $\asym_p$ descends to cohomology
\begin{equation*}
\asym_p\!: e^{\cW(\cL)/z}\GM(\cW(\cL)) \to \C[\![z]\!].
\end{equation*}
Then, one defines the higher residue pairing $P_{\cW(\cL)}\!: \GM(\cW(\cL)) \times \GM(\cW(\cL)) \to \C[\![z]\!]$ as
\begin{equation*}
P_{\cW(\cL)}(\Omega_1, \Omega_2) = \sum_p \overline{\asym_p\left(e^{\cW(\cL)/z} \Omega_1\right)} \asym_p\left(e^{\cW(\cL)/z} \Omega_2 \right),
\end{equation*}
where the sum is over the critical points and we set 
\begin{equation*}
\overline{\asym_p\left(e^{\cW(\cL)/z} \Omega_1\right)} = \asym_p\left(e^{\cW(\cL)/z} \Omega_1\right)\big|_{z\mapsto -z}.
\end{equation*}
By construction, there is a formal power series expansion
\begin{equation*}
P_{\cW(\cL)}(\Omega_1, \Omega_2) = \sum_{k \geq 0} K^{(k)}_{\cW(\cL)}(\Omega_1, \Omega_2) z^k, 
\end{equation*}
such that $K^{(0)}_{\cW(\cL)}(\Omega_1, \Omega_2)$ is the classical Grothendieck residue pairing (see e.g. \cite{GriffithsHarris}, Chapter 5, Section 1), 
\begin{equation*}
K^{(0)}_{\cW(\cL)}(\Omega_1, \Omega_2) = \res_{\cW(\cL)}(\Omega_1, \Omega_2),
\end{equation*}
while $K^{(k)}_{\cW(\cL)}(\Omega_1, \Omega_2)$ are also called higher residue pairings.
\begin{rmk} If $\cW(\cL)$ does not have isolated, nondegenerate critical points, the strategy of \cite{CoatesCortiIritani_hodge} is to pass to a suitable perturbation, by turning on certain equivariant parameters.
\end{rmk}
\subsection{Test configurations, DF invariant and residue pairing}\label{DFResSec}
Let $(\cX, \cL) \to \PP^1$ denote a compactified toric test configuration for a polarised toric manifold $(X, L)$, with smooth total space $\cX$ (see e.g. \cite{Legendre_localised}, Section 4.1). It is known that $K$-semistability can be checked using such test configurations (see \cite{Legendre_localised}, Remark 4.1). 

Since $\cX$ is a smooth projective toric variety, according to our discussion above, it admits a LG mirror family, which we denote by $\cW\!: \widehat{\cY}_{\cX} \to \C$, with respect to a degenerating torus fibration $\widehat{\cY}_{\cX} \to \widehat{\cM}_{\cX}$. Choosing a trivialisation of the torus fibration $\cY_{\cX} \to \cM_{\cX}$ we can write $\cW$ as a Laurent polynomial depending on parameters, 
\begin{equation*}
\cW(x; y) = \sum^{m'}_{i = 1} y_i Q^{\lambda(b'_i)} x^{b'_i} + \sum_{ \kk' \in G_{\cX}} y_{\kk'} Q^{\lambda(\kk')} x^{\kk'}, 
\end{equation*}
where $b'_1, \ldots, b'_{m'}$ are generators of the rays of the fan of $\cX$ in its fan lattice $N_{\cX}$, $Q$ is the Novikov variable, $\lambda(b'_i), \lambda(\kk') \in H_2(\cX, \Q)$ are curve classes, $x \in \Hom(N_{\cX}, \C^*)$ is a torus coordinate, $y_i, y_{\kk'}$ are deformation parameters, and $G_{\cX} \subset N_{\cX}$ is a suitable finite subset. 

We know that, at least for testing $K$-semistability, it is not restrictive to assume that $\cX$ (which is smooth) also admits a toric morphism $\cX \to X \times \PP^1$: this can be seen by replacing $\cX$ by a roof resolving the toric birational morphism $X \times \PP^1 \dashrightarrow \cX$ (see e.g. \cite{Dervan_crit_metrics}, Section 2.1.1). Thus, we have a natural inclusion $H_2(X, \Q) \subset H_2(\cX, \Q)$, a canonical splitting $N_{\cX} \cong N \times \Z$, and the fan of $\cX$ is obtained a a refinement of the fan of $X \times \PP^1$, compatible with the toric morphism $\cX \to X \times \PP^1$. 

This means that we can regard $\cW(x; y)$ naturally as a deformation of $W(x, y)$, 
\begin{equation*}
\cW(x; y) = W(x; y) + \cW'(x; y),
\end{equation*}
where (abusing notation slightly) we still write $x$ for the element of $\Hom(N , \C^*)$ which is obtained by restriction to $N \times \{0\} \subset N \times \Z \cong N_{\cX}$. 
\begin{exm} Suppose $(\cX, \cL_r)$ is given by the degeneration to the normal cone of a point in $(\PP^1, \olo(1))$, with parameter $r \in (0, 1)$. The total space $\cX$ is Fano, given by the del Pezzo $\Bl_{p\times\{0\}} \PP^1\times\PP^1$, and the polarisation $\cL_r : = \olo(1,1) - r E$ is a multiple of $-K_{\cX}$ iff $r = \frac{1}{2}$. Then, the mirror LG family is given by
\begin{equation*}
\cW(\cL_r) = x + \frac{e^{-2\pi}}{x} + \frac{e^{-2\pi}}{x'} + x'(1 + e^{2\pi r} x),
\end{equation*}
where $x + \frac{e^{-2\pi}}{x}$ is the LG potential of $(\PP^1, \olo(1))$, and $x + \frac{e^{-2\pi}}{x} + x' + \frac{e^{-2\pi}}{x'}$ is the LG potential of the trivial test configuration $(\PP^1 \times \PP^1, \olo(1,1))$. Since we will work in the large volume limit, it is also useful to note
\begin{equation*}
\cW(k\cL_r) = x + \frac{e^{-2\pi k}}{x} + \frac{e^{-2\pi k}}{x'} + x'(1 + e^{2\pi r k} x).
\end{equation*}
\end{exm}
\begin{exm} Suppose $\cX$ is isomorphic to the Hirzebruch surface $\PP(\olo\oplus\olo(1)) \cong \Bl_p \PP^2$ endowed with the polarisation $\cL_r = H - r E$ for $r \in (0,1)$. Then $(\cX, \cL_r) \to \PP^1$ is a compactified toric product test configuration for $\PP^1$ (i.e. it is induced by a holomorphic vector field on the fibre $\PP^1$). The total space is del Pezzo and $\cL_r$ is a multiple of the anticanonical iff $r = \frac{1}{3}$ in which case $\cL_r = - \frac{1}{3}K_{\cX}$. A LG potential is given by
\begin{equation*}
\cW(k \cL_r) = \frac{e^{-2\pi k}}{x x'} + x + x' + e^{2\pi k r} x x'.
\end{equation*}
Note that this is not presented as an iterated blowup of $\PP^1 \times \PP^1$ (for that we would need to pass to a blowup of $\cX$).     
\end{exm}
\begin{exm}\label{SlopeUnstableSurfaceExmPotentials} Let $X = \Bl_p \PP^2$ and $(X, L) = (X, -K_X)$. With respect to suitable torus coordinates, the LG potential is given by
\begin{equation*}
W_{X}(-K_X) = x + x' + \frac{1}{x x'} + x x'.
\end{equation*}  

Consider a compactified toric test configuration $(\cX, \cL)$ for $(X, L)$ such that $\cX$ is given by the degeneration to the normal cone of the exceptional divisor $E$, namely $\cX = \Bl_{E \times \{0\}} X \times \PP^1$. Then $\cX$ is a toric Fano threefold of rank $4$, the variety 4-11 in the classification. A standard presentation of $\cX$ is given by the face fan of the reflexive polytope 82 in the Kreuzer-Skarke list \cite{KreuzSkarke},   
\begin{align*}
\cP = \operatorname{conv}(\{\left(1,\,0,\,0\right), \left(0,\,1,\,0\right), \left(0,\,-1,\,0\right), \left(-1,\,0,\,0\right), \left(0,\,0,\,1\right), \left(1,\,1,\,0\right), \left(1,\,0,\,-1\right)\}).
\end{align*}
A corresponding LG potential for $c_1(\cX)$, with respect to the variables of this lattice, is given by
\begin{align*}
& \cW_{\cX}(-K_{\cX}) = X Y + X + Y + Z + \frac{X}{Z} + \frac{1}{Y} + \frac{1}{X} \\
&= \frac{1}{Y} + \left(\frac{X}{Z} + Z + X + \frac{1}{X}\right) + Y + Y X.
\end{align*}
Thus, after the change of variables $x = \frac{X}{Z},\,x' = Z,\,z = Y$, we have
\begin{align*}
\cW_{\cX}(-K_{\cX}) = \frac{1}{z} + \left(x + x' + x x' + \frac{1}{x x'}\right) + z (1+ x x'), 
\end{align*}
giving the required presentation of $\cW_{\cX}(-K_{\cX})$ as a deformation of $W_{X}(-K_X)$ or more precisely of
\begin{equation*}
\cW_{X \times \PP^1}(-K_{X \times \PP^1}) =  \frac{1}{z} + \left(x + x' + x x' + \frac{1}{x x'}\right) + z.
\end{equation*}
\end{exm}

Returning to the general case, let us now consider the Donaldson-Futaki invariant of $(\cX, \cL)$. Recall that, following Odaka \cite{Odaka_blowup} and Wang \cite{Wang_GIT} (see also \cite{Legendre_localised} for the general K\"ahler case), we can write this as a Poincar\'e pairing
\begin{align}\label{OdakaWangFormula}
& \nonumber \DF(\cX, \cL) = \int_{\cX} (c_1(\cL))^n \cup\left(\frac{n c}{n+1}c_1(\cL) + c_1(K_{\cX / \PP^1})\right)\\
& = \left((c_1(\cL))^n, \frac{n c}{n+1}c_1(\cL) + c_1(K_{\cX / \PP^1})\right), 
\end{align}
where
\begin{equation*}
c = \frac{c_1(X) \cup c^{n-1}_1(L)}{c^n_1(L)}.
\end{equation*}
Let us set     
\begin{align}\label{alphaForm}
\nonumber &\tilde{\alpha}_{\cX}(\cL) = \Theta^{-1}_{\cL} \left((c_1(\cL))^n\right) \in \opH^n(\Omega^{\bullet}, z d+ d \cW_{\cX}(\cL)\wedge),\\
&\alpha_{\cX}(\cL) = \Theta^{-1}_{\cL} \left((c_1(\cL))^n\right)\big|_{z=0} \in \opH^n(\Omega^{\bullet}, d \cW_{\cX}(\cL)\wedge), 
\end{align}
and similarly
\begin{align}\label{betaForm}
\nonumber &\tilde{\beta}_{\cX}(\cL) = \Theta^{-1}_{\cL} \left(\frac{n c}{n+1}c_1(\cL) + c_1(K_{\cX / \PP^1})\right) \in \opH^n(\Omega^{\bullet}, z d+ d \cW_{\cX}(\cL)\wedge),\\
&\beta_{\cX}(\cL) =\Theta^{-1}_{\cL} \left(\frac{n c}{n+1}c_1(\cL) + c_1(K_{\cX / \PP^1})\right)\big|_{z=0} \in \opH^n(\Omega^{\bullet}, d \cW_{\cX}(\cL)\wedge).
\end{align}
\begin{lemma}[\cite{Stoppa_LargeComplex}, Section 2] We have 
\begin{align*}
&\DF(\cX, \cL) = P_{\cW_{\cX}(\cL)}\left( \tilde{\alpha}_{\cX}(\cL), \tilde{\beta}_{\cX}(\cL)\right) = K^{(0)}_{\cW_{\cX}(\cL)}(\alpha_{\cX}(\cL), \beta_{\cX}(\cL))\\
& = \res_{\cW(\cL)}(\alpha_{\cX}(\cL), \beta_{\cX}(\cL)).
\end{align*}
\end{lemma}
\begin{proof} Applying the Hodge-theoretic mirror theorem discussed in Section \ref{DFResSec} to $\cX$ shows that we have in particular
\begin{align*}
& \DF(\cX, \cL) = \left((c_1(\cL))^n, \frac{n c}{n+1}c_1(\cL) + c_1(K_{\cX / \PP^1})\right) = P_{\cW(\cL)}\left( \tilde{\alpha}_{\cX}(\cL), \tilde{\beta}_{\cX}(\cL)\right) \in \C[\![z]\!],
\end{align*}
using the linear mirror isomorphism $\Theta_{\cL}\!: \GM(\cW(\cL)) \xrightarrow{\,\,\cong\,\,} H^*(\cX, \C) \otimes \C[\![z]\!]$ and higher residue pairing with respect to the LG potential $\cW_{\cX}(\cL)$ corresponding to $\cL$. 

On the other hand, $\DF(\cX, \cL)$ does not depend on the formal variable $z$, and so must equal the specialisation of the higher residue pairing appearing on the right hand side at $z = 0$. By our discussion in Section \ref{HigherResSec}, this is indeed well defined and equal to the classical Grothendieck residue pairing, defined on the top cohomology of the complex $\left(\Omega^{\bullet}, d \cW(\cL)\wedge\right)$, so we have
\begin{align}\label{DFGrothRes}
\DF(\cX, \cL) = K^{(0)}_{\cW_{\cX}(\cL)}(\alpha_{\cX}(\cL), \beta_{\cX}(\cL)) = \res_{\cW(\cL)}(\alpha_{\cX}(\cL), \beta_{\cX}(\cL)).
\end{align}
\end{proof}

\begin{rmk} The reference \cite{CorreaResidue} shows that the usual Atiyah-Bott localisation formula for the classical Futaki character can also be expressed in terms of certain Grothendieck residue pairings on $X$ (not on the mirror!).  
\end{rmk}
From the viewpoint of algebro-geometric stability, it would be desirable to express the quantities $P_{\cW_{\cX}(\cL)}\left( \tilde{\alpha}_{\cX}(\cL), \tilde{\beta}_{\cX}(\cL)\right)$, $K^{(0)}_{\cW_{\cX}(\cL)}(\alpha_{\cX}(\cL), \beta_{\cX}(\cL))$ as global intersection pairings on $\Hom(N_{\cX}, \C^*)$ or on a suitable compactification. In fact this problem has been studied in different contexts, a priori not related to stability. We will discuss one approach in the next Section.
\subsection{Localisation for the DF invariant}
Let $(\cX^{\vee}, [\eta] , [\xi])$ denote any regular compactified toric test configuration for $X^{\vee} := \cX^{\vee}_1$, allowing a twist $[\xi]$ (as in the Introduction). The Atiyah-Bott localisation formula for the Donaldson-Futaki invariant is discussed by Legendre \cite{Legendre_localised} (the exposition there generalises immediately to allow a twist). In the case of isolated fixed points, it is given by
\begin{align}\label{DFAtiyahBottLocalisation}
&\nonumber\DF( \cX^{\vee}, [\eta] , [\xi])\\
&\nonumber = \sum_{p \in Z(v)}\left(\frac{n c^{\vee}_{\eta, \xi}}{n+1}\frac{(-h_{\eta}(p))^{n+1}}{e(T_p)(v)} - \frac{(\sum^{n+1}_{i=0} w_i(p) +  h_{\xi}(p))(-h(p))^n}{e(T_p)(v)}+ \frac{(-h_{\eta}(p))^n}{e(T_p)(v)}\right)\\
&= \sum_{p \in Z(v)} \frac{(-h_{\eta}(p))^{n}\big(-\frac{nc^{\vee}_{\eta, \xi}}{n+1} h_{\eta}(p)-\sum^{n+1}_{i=0} w_i(p) +1 -  h_{\xi}(p)\big)}{e(T_p)(v)},  
\end{align}
where 
\begin{equation*}
c^{\vee}_{\eta, \xi} := \frac{(c_1(X^{\vee}) - \xi|_{X^{\vee}})\cup (\eta|_{X^{\vee}})^{n-1}}{(\eta|_{X^{\vee}})^n},
\end{equation*}
$h_{\eta}$, $h_{\xi}$ denote the Hamiltonians for the infinitesimal generator $v$ of the structure $\C^*$-action with respect to $\eta$, $\xi$, and the sum is over zeroes $Z(v) = Z(v_i)$ lying in the central fibre $\cX^{\vee}_0$. 
\section{Cohomology intersection form}\label{LogBackSec} 
In the context of logarithmic connections on a smooth projective variety, under certain assumptions, the \emph{cohomology intersection form} provides a global version of the stationary phase expansion for higher residue pairings. Here we summarise a few basic facts which we need for our applications.

Matsumoto \cite{Matsumoto_intersection} studies this intersection form in the classical case of generic hyperplane arrangements in projective space (however he does not discuss the stationary phase formula).  We will follow a recent reference, including the stationary phase formula, due to Matsubara-Heo \cite{Matsubara_cohomology}. We also point out the work of Sabbah \cite{Sabbah_Duality}, and the $L^2$ approach developed by Li and Wen \cite{SiLi_Hodge}. 
 
As in \cite{Matsubara_cohomology}, Section 2.1, we consider a smooth projective variety $Y$ of dimension $n$, endowed with a simple normal crossing divisor $D = \sum^N_j D_j$. Let $E \to Y$ denote a holomorphic vector bundle endowed with a meromorphic integrable connection $\nabla$ with logarithmic poles along $D$ (see e.g. \cite{hotta}, Chapter 5, Section 5.2 for the general theory). Writing $E_+ := E$, $E_- := E^{\vee}$ and $\Omega^p_{\log} := \Omega^p_{Y}(\log D)$ for the sheaf of logarithmic $p$-forms, we set 
\begin{equation*}
\Omega^p_{+} := \Omega^p_{\log} \otimes E_+,\,\Omega^p_{-} := \Omega^p_{\log} \otimes E_-, 
\end{equation*}
and so obtain complexes
\begin{equation*}
\left(\Omega^{\bullet}_{+}, \nabla_+ := \nabla\right),\,\left(\Omega^{\bullet}_{-}, \nabla_- := \nabla^{\vee}\right), 
\end{equation*}
with hypercohomologies 
\begin{equation*}
\opH^p_{\pm} := \HH^p\left(Y; \left(\Omega^{\bullet}_{\pm}, \nabla_{\pm}\right)\right).
\end{equation*}
Under a certain generic condition on $\nabla$, the cohomology intersection form is defined as a bilinear pairing
\begin{equation*}
\bra - , - \ket_{ch}\!: \opH^n_- \otimes_{\C} \opH^n_+ \to \C. 
\end{equation*}
Namely, introduce the conditions 
\begin{equation*}
(!)_+: \Spec(\res_i(\nabla)) \cap \Z_{\leq 0} = \emptyset,\,(!)_-: \Spec(\res_i(\nabla)) \cap \Z_{\geq 0} = \emptyset,  
\end{equation*}
(see Remark \ref{GenericMotivationRmk} for motivation), where $\res_i(\nabla)$ is the endomorphism given by the residue of $\nabla$ along $D_i$, with eigenvalues $\Spec(\res_i(\nabla)) \subset \C$. 
\begin{lemma}[see e.g. \cite{Matsubara_cohomology}, Section 2.1] The condition $(!)_{\pm}$ implies that for all $p$ there is a canonical isomorphism
\begin{equation*}
\reg_{\pm}\!: \opH^p_{\pm} \xrightarrow{\,\cong\,} \opH^p_{\pm}(-D) := \HH^p\left(Y; \left(\Omega^{\bullet}_{\pm}(-D); \nabla_{\pm}\right)\right),
\end{equation*}
where we use the fact that $\nabla_{\pm}$ also induces a connection on $E_{\pm}(j D)$ with logarithmic poles along $D$ for all $j \in \Z$. 
\end{lemma}
\begin{rmk}\label{GenericMotivationRmk} This is a complex-analytic analogue of the fact, used in the classical reference \cite{Matsumoto_intersection}, that for complements of hyperplane arrangements in projective space, at least under the stronger condition that both $(!)_{+}$ and $(!)_{-}$ hold, there is a natural isomorphism
\begin{equation*}
H^k_c(Y, \nabla) \xrightarrow{\,\,\cong\,\,} H^k(Y, \nabla), 
\end{equation*}
using twisted de Rham cohomology groups defined using compactly supported smooth forms and arbitrary smooth forms, respectively. 

The necessity of one of the conditions $(!)_{\pm}$ can be checked e.g. by regarding $d - \frac{dz}{z}$ as a meromorphic connection $\nabla$ on the trivial line bundle on $Y = \PP^1$, with logarithmic poles at $D = 0 + \infty$, with $z$ denoting a standard coordinate on $\C^*$. The spaces $\opH^0_{+}$, $\opH^0_{+}(-D)$ are not isomorphic, as witnessed by the flat section given by the polynomial $z$.    
\end{rmk}
\begin{rmk} In fact if both $(!)_{+}$ and $(!)_{-}$ hold then we must have 
$\opH^p_{\pm} = 0,\,p \neq n$.
\end{rmk}
By the usual Dolbeault argument, the canonical morphisms of complexes 
\begin{equation*}
\left(\Omega^{\bullet}_{\pm}(j D), \nabla_{\pm}\right) \to \left(\mathcal{E}^{\bullet}_{\pm}(j D), \nabla_{\pm} + \delbar\right), 
\end{equation*}  
the latter involving the sheaves of \emph{smooth} forms $\mathcal{E}^{p, q}_{\pm} := \Omega^p_{\pm} \otimes_{\olo_Y} \mathcal{E}^{0,q}_Y$, are quasi-isomor\-phisms. This means in particular that we can choose smooth $n$-forms $\omega_{-}$ on $Y \setminus D$ and $\omega_+$ on $Y$ representing given classes $[\omega_-] \in \opH^n_-$, respectively $[\omega_+] \in \opH^n_+(-D)$, under the Dolbeault isomorphism. So we have a well defined pairing
\begin{align*}
& \bra - , - \ket_+\!: \opH^n_-\otimes_{\C} \opH^n_+(-D) \to \C,\, \bra [\omega_-], [\omega_+] \ket_+ := \left(\frac{1}{2\pi \ii}\right)^n \int_Y \omega_- \wedge \omega_+ 
\end{align*}
(where we are also taking the trace of the endomorphism part). 

We can summarise our discussion as follows.
\begin{prop}[see e.g. \cite{Matsubara_cohomology}, Section 2.1] Suppose that the condition $(!)_{+}$ holds. Then there is a well defined \emph{cohomology intersection form} 
\begin{align*}
\bra - , - \ket_{ch}\!: \opH^n_-\otimes_{\C} \opH^n_+ \to \C  
\end{align*}
given by
\begin{align*}
\bra [\omega_-], [\omega_+] \ket_+ := (2\pi \ii)^n\bra [\omega_-], \reg_+ [\omega_+] \ket_+  
\end{align*}
(as we noted, replacing $[\omega_+]$ by $\reg_+ [\omega_+]$ is analogous to choosing a compactly supported representative). In particular, we obtain a cohomology intersection form
\begin{align*}
& \bra - , - \ket_{ch}\!: \opH^0(Y, \Omega^n_-) \otimes_{\C} \opH^0(Y, \Omega^n_+) \to \C 
\end{align*}
by using the natural maps $\opH^0(Y, \Omega^n_{\pm}) \to \opH^n_{\pm}$. Similar statements hold when the condition $(!)_{+}$ is replaced with $(!)_{-}$.
\end{prop}
\subsection{Residue theorem}\label{ResThmSubSec}
In the classical case of generic hyperplane arrangements in projective space, Matsumoto \cite{Matsumoto_intersection}, Theorem 2.1, proves a \emph{residue theorem} for the cohomology intersection form, that is, a formula computing $\bra [\omega_-] , [\omega_+] \ket_{ch}$ in terms of contributions from strata of the hyperplane arrangement. 

A general version is proved in \cite{Matsubara_cohomology}. Assume the generic condition $(!)_+$ for definiteness (this could be replaced with $(!)_-$). Fix logarithmic $n$-forms $\omega_{\pm} \in \opH^0(Y, \Omega^n_{\pm})$. For a fixed ordered multi-index $P_n$ of length $n$ (i.e. an ordered partition with $n$ components), write 
\begin{equation*}
D(P_n) := \cap_{i} D_{P_n(i)}
\end{equation*}
for the corresponding $0$-dimensional stratum of the boundary divisor $D$. There is a natural notion of restriction $\res_{P_n}(\omega)$, i.e. factoring out the simple poles in the given order. Namely, decomposing $\omega = \frac{dx_i}{x_i} \wedge \omega' + \omega''$ near $D_i = \{x_i = 0\}$, we set 
\begin{equation*}
\res_i(\omega) := \omega'|_{D_i}
\end{equation*}
and extend this operation by 
\begin{equation*}
\res_{P_n}(\omega) := \res_{P_n(n)}\circ \cdots \circ \res_{P_n(1)}(\omega). 
\end{equation*}

The residue of the connection $\res_{P_n}(\nabla)$ is defined in the same way (the ordering induced by $P_n$ is irrelevant in this case, by integrability). 
\begin{rmk}
We note that the (potentially confusing) notation $\res_{P_n}(\omega)$ for the \emph{restriction of a meromorphic form}, respectively $\res_{P_n}(\nabla)$ for the \emph{residue of a meromorphic connections} is chosen to agree with \cite{Matsubara_cohomology}, Section 2.1; this is justified by the similarity between the two operations (since $\res_{P_n}(\omega)$ may be regarded as a form-valued residue).
\end{rmk} 
Define
\begin{equation*}
\langle \res_{P_n}(\omega_+)\big| \res_{P_n}(\nabla)^{-1}\big| \res_{P_n}(\omega_-)\rangle := \sum_{z \in D(P_n)} \langle \res_{P_n}(\nabla)^{-1}_z \res_{P_n}(\omega_+), \res_{P_n}(\omega_-)_z\rangle,
\end{equation*}
where the right hand side uses the duality bewteen $E_z$, $E^{\vee}_z$. Then we have
\begin{thm}[Residue theorem, \cite{Matsubara_cohomology} Theorem 2.2]\label{MainResThm} Suppose that the condition $(!)_{+}$ holds. Then we have
\begin{equation*}
\bra \omega_-, \omega_+\ket_{ch} = (-2\pi \ii)^n \sum_{P_n : D(P_n) \neq \emptyset} \langle \res_{P_n}(\omega_+)\big| \res_{P_n}(\nabla)^{-1}\big| \res_{P_n}(\omega_-)\rangle.
\end{equation*}
\end{thm}
\subsection{Stationary phase formula}\label{StationaryPhaseSubSec}
Suppose that the connection $\nabla$ has rank $1$, i.e. $E \to Y$ is a line bundle, and that the complement $U:= Y\setminus D$ is affine. There is corresponding $1$-parameter family of deformed connections $\nabla^z$ such that its restriction to $U$ is given by $\nabla^z = z d + \alpha\wedge$ for a fixed connection form $\alpha$. We continue to assume the generic condition $(!)_+$ for definiteness (this could be replaced with $(!)_-$)
\begin{definition}[see e.g. \cite{Matsubara_cohomology}, Section 2.2]\label{StatPhasDef} We say that \emph{the stationary phase formula holds for} $(Y, \nabla)$ if the condition $(!)_{+}$ holds and there is a Laurent series expansion for the cohomology intersection pairing with respect to the deformed connections,
\begin{align*}
\bra \omega_-, \omega_+\ket_{ch, \nabla^z} = ( 2\pi \ii z)^n \sum_{k \geq 0} K^{(k)}(\omega_-, \omega_+) z^{k},\,z \to 0,
\end{align*}
such that $K^{(0)}(\omega_-, \omega_+)$ is the Grothendieck residue pairing of the $n$-forms $\omega_-$, $\omega_+$ on $U$.
\end{definition}
\begin{lemma}[\cite{Matsubara_cohomology}, Corollary 1.4]\label{HomogeneityLem} If the stationary phase formula holds for $(Y, \nabla)$, then in fact we must have
\begin{align*}
\bra \omega_-, \omega_+\ket_{ch, \nabla^z} = ( 2\pi \ii)^n K^{(0)}(\omega_-, \omega_+),
\end{align*}
that is, all higher residues must vanish. In particular, by specialising at $z = 1$, we obtain, for the original cohomology intersection pairing,
\begin{equation*}
\bra \omega_-, \omega_+\ket_{ch} = ( 2\pi \ii)^n K^{(0)}(\omega_-, \omega_+).
\end{equation*}
\end{lemma}
\begin{proof}
The result follows at once by applying Theorem \ref{MainResThm} to the deformed connections $\nabla^z$ (rather than the undeformed $\nabla$), since this implies that $\bra \omega_-, \omega_+\ket_{ch, \nabla^z}$ is homogeneous of degree $n$ with respect to $z$. 
\end{proof}

Thus, if we are only interested in a Grothendieck residue $K^{(0)}(\omega_-, \omega_+)$ (which is the case for our expression \eqref{DFGrothRes} for the Donaldson-Futaki invariant), we could hope that, after a suitable extension to logarithmic forms, a stationary phase formula holds, allowing the application of the residue theorem \ref{MainResThm}. This is the course we will follow in the next Sections.

The stationary phase formula for logarithmic connections is proved in \cite{Matsubara_cohomology}, Theorem 1.1, under some structure conditions on the connection form $\alpha$ (spelled out in \cite{Matsubara_cohomology}, Section 2.2; in particular, the line bundle $E \to Y$ should be trivial). We note that this is closely related to work of Huh \cite{huh_likelihood}. 

Rather than recalling these more general structure conditions, here we will only state a \emph{sufficient condition}, of a topological nature, under which they hold automatically; according to \cite{Matsubara_cohomology}, Remark 2.6 (5), this follows from results of Huh \cite{huh_likelihood}.

Thus, as in \cite{Matsubara_cohomology}, Section 2.2, we consider a logarithmic connection $\nabla = d + F \wedge $ on $U = Y \setminus D$, where 
\begin{equation*}
F = \sum^m_{i = 1} \alpha_i \log f_i,
\end{equation*} 
and the $f_i$ are regular functions on $U$. Suppose that 
\begin{itemize}
\item[$(M)$] the critical locus $\Crit(F)$ is discrete and the sum of Milnor numbers 
\begin{equation*}
m_p := \dim_{\C} \olo_{U, p}/(\del_{x_1} F, \ldots, \del_{x_n} F)
\end{equation*}
satisfies 
\begin{equation*}
\sum_{p \in \Crit(F)} m_p = (-1)^n \chi(U).
\end{equation*}
\end{itemize}

\begin{lemma}[see \cite{Matsubara_cohomology}, Theorem 2.5, Remark 2.6 (5), and Corollary 2.7]\label{StatPhasLem} Under the condition $(M)$, the stationary phase formula holds in the sense of Definition \ref{StatPhasDef}, so that for $\omega_{\mp} \in \opH^0(U, \Omega^n_U)$ we have
\begin{align*} 
\left( \frac{1}{2\pi \ii}\right)^n\bra \omega_-, \omega_+\ket_{ch} = K^{(0)}(\omega_-, \omega_+).
\end{align*}
Thus, by the classical expression for the residue pairing $K^{(0)}$ (see e.g. \cite{GriffithsHarris}, Chapter 5, Section 1), we have
\begin{align*} 
\left( \frac{1}{2\pi \ii}\right)^n\bra \omega_-, \omega_+\ket_{ch} = \left( \frac{1}{2\pi \ii}\right)^n \sum_{p \in \Crit(F)} \int_{\Gamma_p} \frac{\frac{\omega_+}{dx}\frac{\omega_-}{dx}}{\del_{x_1} F \cdots \del_{x_n} F } dx,
\end{align*}
where $(x_1, \ldots, x_n)$ denote local coordinates in a neighbourhood of a critical point $p$, we set $dx = dx_1 \wedge \cdots \wedge dx_n$, and the integration cycle is given by $\Gamma_p = \{|\del_{x_i} F| = \varepsilon, i = 1,\ldots, n\}$ for sufficiently small $\varepsilon > 0$, oriented so that $d\arg |\del_{x_1} F| \wedge \cdots \wedge d\arg |\del_{x_n} F| > 0$.  
\end{lemma}

\section{Application of $\bra - , - \ket_{ch}$ to toric LG models}\label{MainSec}
\subsection{DF invariants and $\bra - , - \ket_{ch}$ on the mirror}\label{DFchSubSec} Let $\cW := \cW_{\cX}(\cL)$ denote a Laurent polynomial on the affine torus $\T := (\C^*)^{n+1}$ which is mirror to the polarised toric manifold $(\cX, \cL)$ (in a fixed trivialisation of the mirror family $\widehat{\cY}_{\cX} \to \widehat{\cM}_{\cX}$ on the generic locus). 

We introduce auxiliary logarithmic connections on $\T \setminus V(\cW)$ given by
\begin{align*}
\nabla_{\cW} = d + d\log\cW \wedge,\, \nabla^z_{\cW} := z d + d\log\cW \wedge. 
\end{align*}
In the light of our discussion in the previous Sections, our aim will be to obtain a global expression for the Donaldson-Futaki invariant as a cohomology intersection number,
\begin{equation*}
\DF(\cX, \cL) = \left(\frac{1}{2\pi \ii z} \right)^{n+1}\bra \tilde{\omega}_-, \tilde{\omega}_+ \ket_{ch, \nabla^z_{\cW}},
\end{equation*}
where $[\tilde{\omega}_{\pm}] \in \opH^{n+1}_{\pm}$ are now chosen as classes of logarithmic $(n+1)$-forms extending the classes 
\begin{align*}
 \cW^{-n}\overline{\tilde{\alpha}_{\cX}(\cL)},\,\cW^{-1}\tilde{\beta}_{\cX}(\cL)
\end{align*}
and we set $\overline{\tilde{\alpha}_{\cX}(\cL)} := \tilde{\alpha}_{\cX}(\cL)|_{z\mapsto - z}$.

Several conditions are required for this. Firstly, it will be convenient to allow the prospective logarithmic form $\cW^{-n}\overline{\tilde{\alpha}_{\cX}(\cL)}$ to be replaced by $\cW^{-1}_{n, \eps}\overline{\tilde{\alpha}_{\cX}(\cL)}$, where
\begin{equation*}
\cW_{n, \eps} := \prod^{n}_{i = 1} \cW^{(i)}_{\eps},
\end{equation*}
and each $\cW^{(i)}_{\eps}$ is a generic analytic perturbation of $\cW$ parametrised by $\eps$ in the unit disc. Accordingly, we will work on $\cT := \T \setminus (V(\cW) \cup V(\cW_{n, \eps}))$. This will simplify the choice of a suitable compactification $\barcT$. This compactification should satisfy some key properties.  
\begin{definition}\label{CompactPropertiesDef} Let us introduce conditions for a compactification $\barcT$ of $\cT$ given by:
\begin{enumerate}
\item[$(i)$] $\barcT$ is smooth, it induces a compactification of the structure torus $\T$, and the complement $\barcT \setminus \cT$ is a simple normal crossings divisor $\cD \subset \T$;     
\item[$(ii)$] the rank $1$ holomorphic integrable connection $\nabla^z_{\cW}$ extends to a meromorphic flat connection on the trivial line bundle over $\barcT$ with logarithmic poles along $\cD$;
\item[$(iii)$] the classes of holomorphic $n+1$-forms $\cW^{-1}_{n, \eps}\overline{\tilde{\alpha}_{\cX}(\cL)}, \cW^{-1}\tilde{\beta}_{\cX}(\cL)$ extend to classes of logarithmic $n+1$-forms $[\tilde{\omega}_{\pm}] \in \opH^{n+1}_{\pm}$ on $\barcT$. We will write $\omega_{\pm} := \tilde{\omega}_{\pm}|_{z=0}$ for their specialisation;
\item[$(iv)$] the stationary phase formula holds for the corresponding cohomology intersection form $\bra - , -\ket_{ch, \nabla^z_{\cW}}$, computed on $\barcT$ with respect to $\nabla^z_{\cW}$, in the sense of Definition \ref{StatPhasDef}. 
\end{enumerate}
\end{definition}
\begin{lemma}\label{FutakiToCohoIF} Suppose $\barcT$ is a compactification of $\cT$ satisfying the conditions of Definition \ref{CompactPropertiesDef}. Then, we have
\begin{align*}
\DF(\cX, \cL) = \left( \frac{1}{2\pi \ii z} \right)^{n+1}\bra \tilde{\omega}_-, \tilde{\omega}_+ \ket_{ch, \nabla^z_{\cW}} + O(\eps) = \left( \frac{1}{2\pi \ii}\right)^{n+1} \bra \omega_-,\omega_+ \ket_{ch} + O(\eps).
\end{align*}
\end{lemma}
\begin{proof} Let us set $dx = \prod^{n+1}_{i=1} dx_i$, $\frac{dx}{x} = \prod^{n+1}_{i=1} \frac{dx_i}{x_i}$. By the conditions of Definition \ref{CompactPropertiesDef}, we may apply Lemma \ref{HomogeneityLem}, so we have 
\begin{align*}
&\left( \frac{1}{2\pi \ii z} \right)^{n+1}\bra \tilde{\omega}_-, \tilde{\omega}_+ \ket_{ch, \nabla^z_{\cW}} = \left( \frac{1}{2\pi \ii}\right)^{n+1} \bra \omega_-,\omega_+ \ket_{ch}\\
&= \res_{\log\cW}(\omega_-,\omega_+)\\
&= \left( \frac{1}{2\pi \ii}\right)^{n+1} \sum_{p \in \Crit(\log\cW)} \int_{\Gamma_p} \frac{\frac{\omega_+}{dx}\frac{\omega_-}{dx}}{\del_{x_1} (\log \cW) \cdots \del_{x_{n+1}} (\log \cW) } dx. 
\end{align*}
It follows from the continuity of the Grothendieck residue that
\begin{align*}
&\lim_{\eps \to 0}\left( \frac{1}{2\pi \ii z} \right)^{n+1}\bra \tilde{\omega}_-, \tilde{\omega}_+ \ket_{ch, \nabla^z_{\cW}} = \left( \frac{1}{2\pi \ii}\right)^{n+1} \sum_{p \in \Crit(\cW)} \int_{\Gamma_p} \frac{\frac{\alpha_{\cX}}{dx}\frac{\beta_{\cX}}{dx}}{\del_{x_1} (\cW) \cdots \del_{x_{n+1}} (\cW) } dx\\
&=\left( \frac{1}{2\pi \ii}\right)^{n+1} \sum_{p \in \Crit(\cW)} \int_{\Gamma_p} \frac{\frac{\alpha_{\cX}}{(dx/x)}\frac{\beta_{\cX}}{(dx/x)}}{x_1\del_{x_1} (\cW) \cdots x_{n+1}\del_{x_{n+1}} (\cW) } \frac{dx}{x}\\
&= \res_{\cW}(\alpha_{\cX}(\cL), \beta_{\cX}(\cL))\\
&= \DF(\cX, \cL)
\end{align*}
by \eqref{DFGrothRes}.
\end{proof}

In Proposition \ref{GoodCompactProp}, we will construct a compactification $\barcT$ satisfying conditions $(i)$, $(ii)$, $(iii)$ of Definition \ref{CompactPropertiesDef}, at least nearby the large volume limit, to order $O(k^{-1})$. As a preliminary step we construct an intermediate toric compactification satisfying suitable properties. 
\begin{lemma}\label{IntermediateCompactLem} There exists a (non-unique) smooth toric compactification $\barcT_{\T}$ of the structure torus $\T$, such that $\cW$ extends to a section of a line bundle, and $\nabla^z_{\cW}$ extends to a meromorphic flat connection on the trivial line bundle over $\barcT_{\T}$, with logarithmic poles along $V(\cW)$ and the toric boundary $\cD_{\T}$, possibly away from $\cD_{\T} \cap V(\cW)$ (if this intersection is nonreduced).  
\end{lemma}
\begin{proof} Recall we have a presentation
\begin{equation*}
\cW(\cL)(x; y) = \sum^{m'}_{i = 1} y_i Q^{\lambda(b'_i)} x^{b'_i} + \sum_{ \kk' \in G_{\cX}} y_{\kk'} Q^{\lambda(\kk')} x^{\kk'}, 
\end{equation*} 
introduced in Section \ref{DFResSec}. Note that we can write
\begin{equation*}
\nabla^z_{\cW} = z d + d\log\cW \wedge = z d + \sum^{n+1}_{i=1}\cW^{-1} x_i \del_{x_i} \cW \frac{dx_i}{x_i} \wedge.
\end{equation*}
Choose an intermediate compactification $\barcT'_{\T}$ of $\cT$ given by a smooth \emph{toric} compactification of $\T$, with toric boundary $\cD'_{\T}$, satisfying the property: 
\begin{itemize}
\item[$(T)$] all the monomials $x^{b'_i}$, $i = 1, \ldots, m'$ and  $x^{\kk'}$, $\kk' \in G_{\cX}$ on $\T$ appearing in the LG potential $\cW$ extend to sections of a holomorphic toric line bundle $E \to \barcT'_{\T}$.
\end{itemize}
By the momentum construction of toric varieties, the condition $(T)$ can be achieved by choosing a Delzant momentum polytope containing the exponents $\{b'_i\}$, $i = 1, \ldots, m'$ and  $\{\kk'\}$, $\kk' \in G_{\cX}$. 

Moreover, given a choice of $\barcT'_{\T}$ satisfying $(T)$, by (repeatedly) blowing up torus fixed points, we can pass to an intermediate smooth toric compactification $\barcT_{\T}$, with toric boundary $\cD_{\T}$, satisfying the following condition:
\begin{itemize}
\item[$(T')$] the divisor $V(\cW) \subset \barcT_{\T}$ (which is well-defined by $(T)$) does not contain torus fixed points of $\barcT_{\T}$.  
\end{itemize}
By $(T)$ and $(T')$, this choice of $\barcT_{\T}$ satisfies the properties:
\begin{enumerate}
\item[$(1)$]  by the toric condition, the $1$-forms on $\cT$ given by $d\log(x_i)$, $i =1, \dots, n+1$ extend to meromorphic forms on $\barcT_{\T}$ with a simple pole along $\cD_{\T}$; 
\item[$(2)$] the connection $1$-form $\sum^{n+1}_{i=1}x_i \del_{x_i} \log \cW \frac{dx_i}{x_i}$ extends to a meromorphic $1$-form on the trivial bundle over $\barcT_{\T}$, with logarithmic poles along $\cD_{\T}$ and $V(\cW)$, possibly away from $\cD_{\T} \cap V(\cW)$ (if this intersection is nonreduced). 
\end{enumerate}
This completes the proof of our claims.
\end{proof}
Given its importance for us, we illustrate the construction of $\barcT_{\T}$ in some examples.
\begin{exm} Let $(\cX, k\cL_r)$ be degeneration to the normal cone of a point in $\PP^1$. The convex envelope of the exponents of the monomials
\begin{equation*}
\{x,\,\frac{e^{-2\pi k}}{x},\,x', \frac{e^{-2\pi k}}{x'},\,e^{2\pi k r} x x'\} 
\end{equation*}
is not Delzant (it is the momentum polygon of an orbifold). However, it is contained in the Delzant polygon given by the convex envelope of the exponents of the monomials
\begin{equation*}
\{x,\,\frac{e^{-2\pi k}}{x},\,x', \frac{e^{-2\pi k}}{x'},\,e^{2\pi k r} x x',\, \frac{1}{x x'}\}. 
\end{equation*}
This is the momentum polygon of the toric del Pezzo $S_6$, the blowup of $\PP^2$ at the torus fixed points, with respect to $-K_{S_6}$, thus the monomials of $\cW(k\cL_r)$ correspond to anticanonical sections on $S_6$. With a suitable choice of coordinates for the anticanonical embedding $S_6 \subset \PP[x_0 : \cdots: x_6]$, we have  
\begin{align*}
\cW(k\cL_r) = x_1 + e^{-2\pi k} x_4 +  e^{-2\pi k} x_6 + x_3 + e^{2\pi r k} x_2.
\end{align*}
By standard results in toric geometry, the toric boundary is the union of the smooth rational curves $C_i = \PP[x_i, x_{i+1}]$, using a cyclic index $i = 1, \ldots, 6$, and so the torus fixed points are $p_i := C_i \cap C_{i+1}$. We find that in this case the torus fixed point $p_4$ is contained in the anticanonical divisor $V(\cW(k\cL_r))$, and so an admissible choice of intermediate smooth toric compactification is given by $\barcT_{\T} := \Bl_{p_4} S_6$. 
\end{exm}
\begin{exm} Let $(\cX, \cL_r)$ be the product test configuration given by $(\PP(\olo\oplus\olo(1)) \cong \Bl_p \PP^2, H - r E)$. As in the previous example, the convex envelope of the monomials 
\begin{equation*}
\{ \frac{e^{-2\pi k}}{x x'},\, x,\, x',\, e^{2\pi k r} x x' \} 
\end{equation*} 
is not Delzant, but is contained in the momentum polygon of $(S_6, -K_{S_6})$, and in the same homogeneous coordinates we have
\begin{equation*}
\cW(k\cL_r) =  e^{-2\pi k} x_5 + x_1 + x_3 + e^{2\pi k r} x_2.
\end{equation*}
The torus fixed points $p_3$, $p_5$ are contained in the locus $V(\cW(k\cL_r))$, and an admissible choice of intermediate smooth toric compactification is given by $\barcT_{\T} := \Bl_{p_3, p_5} S_6$.
\end{exm}
\begin{exm} Let $(\cX, -K_{\cX})$ be the toric test configuration for $X = \Bl_p \PP^2$ given by the degeneration to the normal cone of $E$, $\cX = \Bl_{E \times \{0\}} X \times \PP^1$. Then $\cX$ is Fano, given by the face fan of the reflexive polytope  
\begin{align*}
\cP' = \conv(\{(0,0,-1),(1,0,0),(0,1,0),(1,1,0),(-1,-1,0),(0,0,1),(1,1,1)\}).
\end{align*}
The polar dual of $\cP'$ is given by 
\begin{align*}
&(\cP')^{\circ} = \conv(l),\\
&l = \{( 2, -1, -1), ( 2, -1,  1), (-1,  2,  1), (-1,  2, -1),\\
& ( 0, -1,  1), (-1,  0,  1), (-1,  0,  0), ( 0, -1,  0), ( 1, -1, -1), (-1,  1, -1)\}. 
\end{align*}
In order to find a choice of $\barcT_{\T}$ we proceed slightly differently from the previous examples. By the general theory for three-dimensional reflexive polytopes (see \cite{Przyjalkowski_CYLG}), we know that the LG potential $\cW_{\cX}(-K_{\cX})$ corresponds to an anticanonical section on the Gorenstein toric Fano $T_{(\cP')^{\circ}}$. Moreover, $T_{(\cP')^{\circ}}$ admits a (non-unique) crepant resolution induced by a maximal triangulation $\widetilde{T}_{(\cP')^{\circ}}$, and there is a natural anticanonical pencil on $\widetilde{T}_{(\cP')^{\circ}}$ generated by $W_{\cX}(-K_{\cX})$ and the toric boundary. Thus, in this case, we can take 
\begin{equation*}
\barcT_{\T} := \widetilde{T}_{(\cP')^{\circ}}.
\end{equation*}
Note that using, for example, the $GL(3, \Z)$-equivalence $\cP \sim \cP'$ with the polytope $\cP$ in the Kreuzer-Skarke list, one can check that $T_{(\cP')^{\circ}}$ is the Gorenstein canonical Fano threefold with degree $10$ and Picard index $1$ given by the reflexive polytope $\cP^{\circ}$ with index $4185$.
\end{exm}
 
\begin{prop}\label{GoodCompactProp} Replace $\cL$ by a multiple $k \cL$, $k \gg 1$. Then there exists a compactification $\barcT$ of the structure torus $\T$ satisfying the conditions $(i)$, $(ii)$, $(iii)$ of Definition \ref{CompactPropertiesDef}, up to terms of order $O(k^{-1})$. Moreover, there is a canonical embedding of the set of torus fixed points $\cD^{[0]}_{\T} \subset \barcT_{\T}$ into the zero-dimensional stratum $\cD^{[0]}$.
\end{prop}
\begin{proof} We start by recalling some facts concerning the Gauss-Manin connection $\nabla^{\GM}$ acting on $\GM(\cW) = \opH^{n+1}\left(\Omega^{\bullet}_{\widehat{\cY}/\widehat{\cM}}, z d+ d \cW\wedge\right)$. Following \cite{Iritani_survey}, Section 4, we regard the Gauss-Manin connection as a map
\begin{equation*}
\nabla^{\GM}\!: \GM(\cW) \to \frac{1}{z} \GM(\cW) \otimes_{\olo_{\widehat{\cM}_{\cX}}} \Omega^1_{\widehat{\cM}_{\cX}} \oplus \GM(\cW)\frac{dz}{z^2}. 
\end{equation*} 
Since $\cX$ is toric, it can be presented as a GIT quotient 
\begin{equation*}
\cX = \cX_{\xi} := (\C^*)^m /\!\!/_{\xi} K
\end{equation*}
for a torus $K \cong (\C^*)^k$, with Lie algebra $\frk$, contained in the maximal torus $K \hookrightarrow (\C^*)^m$. Here, writing $\tilde{\cD}_1, \ldots, \tilde{\cD}_m \subset \Hom(K, \C^*)$ for the components of the latter embedding, we denote by $\xi \in \frk^{\vee}_{\R}$ a ``K\"ahler" or ``stability" parameter $\xi \in \sum^m_{i = 1}\R_{\geq 0} \tilde{\cD}_i$. Choose a splitting of the dual sequence $1 \to \check{T} \to (\C^*)^m \to \check{K} \to 1$, as well as coordinates $x = (x_1, \ldots, x_{n+1})$ on $\check{T} \cong (\C^*)^{n+1}$ and $q = (q_1,\ldots, q_k)$ on $\check{K} \cong (\C^*)^{k}$. Write $\del_a = q_a \del_{q_a}$, and let $\Omega_0$ denote the standard relative volume form of the family $\widehat{\cY}_{\cX} \to \widehat{\cM}_{\cX}$,
\begin{equation*}
\Omega_0 = \frac{dx_1}{x_1} \wedge \cdots \wedge \frac{dx_{n+1}}{x_{n+1}}.
\end{equation*}
Then, the Gauss-Manin connection acts by
\begin{align}\label{ExplicitGM}
\nabla^{\GM}\big(f \Omega_0\big) &= \sum^k_{a=1} \left(\left(\del_a f + \frac{\del_a \cW}{z}f\right)\Omega_0\right)\frac{dq_a}{q_a}+ \left(\left(z\del_z f - \frac{\cW}{z}f -\frac{n}{2} f\right)\Omega_0\right)\frac{dz}{z}.
\end{align}  

Moreover, according to \cite{Iritani_survey}, Section 5, over the open torus $\check{K} \subset \widehat{\cM}_{\cX}$, $\GM(\cW)$ is generated by $\Omega_0$ as a module over the ring $\olo_{\check{K}}[z]\bra z\del_1, \ldots, z\del_k\ket$, where $\del_a$ acts by the Gauss-Manin connection $\nabla^{\GM}_{\del_a}$. 

We can now turn to achieving properties $(i)$, $(ii)$, $(iii)$ of Definition \ref{CompactPropertiesDef}. The compactification $\barcT$ will be given by a suitable non-toric blowup of the intermediate toric compactification $\barcT_{\T}$ constructed in Lemma \ref{IntermediateCompactLem}. We adopt the notation introduced in the proof of that Lemma. 

Recall that the connection $1$-form $\sum^{n+1}_{i=1}x_i \del_{x_i} \log \cW \frac{dx_i}{x_i}$ extends to a meromorphic $1$-form on the trivial bundle over $\barcT_{\T}$, with logarithmic poles along $\cD_{\T}$ and $V(\cW)$, possibly away from $\cD_{\T} \cap V(\cW)$ (if this intersection is nonreduced).

Moreover, by the toric condition, we know that $\Omega_0$ extends to a meromorphic form on $\barcT_{\T}$ with simple poles along $\cD_{\T}$.

Now we observe that, over $\check{K} \subset \widehat{\cM}_{\cX}$, the $n+1$-form $\tilde{\beta}_{\cX}(\cL)$ can be expressed as a linear combination of $\nabla^{\GM}_{\del_a}(\Omega_0)$, $a = 1, \ldots, k$. By \eqref{ExplicitGM} and the property $(T)$, each $n+1$-form $\nabla^{\GM}_{\del_a}(\Omega_0)$ extends to a meromorphic form on $\barcT$, with values in $E$, with simple poles along $\cD_{\T}$ and $V(\cW)$, possibly away from $\cD \cap V(\cW)$ (if this is nonreduced). Therefore, $\cW^{-1}\tilde{\beta}_{\cX}(\cL)$ extends to a logarithmic $n+1$-forms on $\barcT_{\T}$, with simple poles along $\cD_{\T}$ and $V(\cW)$, away from the base locus of the linear system defined by $\tilde{\beta}_{\cX}(\cL)$ and $\cW$.

However, this argument needs to be refined in order to show that $\cW^{-1}_{n, \eps}\alpha_{\cX}(\cL)$ extends to a logarithmic $n+1$-form on $\barcT_{\T}$, with simple poles along $\cD_{\T}$ and $V(\cW_{n, \eps})$, away from the base locus of the linear system defined by the divisors $(\alpha_{\cX}(\cL))$ and $(\cW_{n, \eps})$, at least after replacing $\cL$ by $k\cL$ for $k \gg 1$, and working modulo terms of order $O(k^{-1})$. 

In order to show this, we consider the module $\GM(\cW(k\cL))|_{z=0}$ given by the top cohomology of the complex $\left(\Omega^{\bullet}, d \cW(k\cL)\wedge\right)$. 

As recalled in Section \ref{ToricLGSec}, choosing the generator $\Theta^{-1}_{k\cL}(1)$ for the $\GM(\cW(k\cL))|_{z=0}$ induces an algebra structure on $\GM(k\cW(\cL))|_{z=0}$ and an isomorphism $\Jac(\cW(k\cL)) \cong (H^*(\cX, \C), *_{[c_1(k\cL)]})$ with the quantum cohomology ring, where we write $*_{[c_1(k\cL)]}$ for the quantum product evaluated at the quantum parameter $Q$, such that
\begin{equation*} 
Q^d = \exp\left(-2 \pi k\int_d c_1(\cL)\right),\,d\in \Lambda_+.
\end{equation*}
Then, we have
\begin{align*}
& k^{-n}\Theta^{-1}_{k\cL} \left((c_1(k\cL))^n\right)\big|_{z=0} = k^{-n}\Theta^{-1}_{k\cL}\left(c_1(k\cL) *_{[c_1(k\cL)]} \cdots *_{[c_1(k\cL)]}c_1(k\cL) + O(k^{-1}) \right)\\
& = \left(k^{-1}\Theta^{-1}_{k\cL}(c_1(k\cL))\right)^n + O(k^{-1}) ,
\end{align*}
since the mirror map intertwines the quantum product with the product on $\Jac(\cW(k\cL))$, and the quantum product approaches the cup product in the large volume limit. Now the same argument used for $\cW^{-1}\tilde{\beta}_{\cX}(\cL)$ shows that $\cW^{-1}_{n, \eps}\left(\Theta^{-1}_{k\cL}(c_1(k\cL))\right)^n$ admits the required logarithmic extension.

Given these extensions, we can construct a simple normal crossing compactification $\barcT$ as an iterated blowup of $\barcT_{\T}$ along $\cD_{\T} \cap V(\cW)$ (if this is nonreduced), and resolving the singularities of the divisors $(\cW(\cL))$, $(\cW_{n, \eps})$, $(\alpha_{\cX}(\cL))$, $(\tilde{\beta}_{\cX}(\cL))$ as well as the base loci $\Bs(\tilde{\beta}_{\cX}(\cL), \cW)$, $\Bs(\alpha_{\cX}(\cL), \cW_{n, \eps})$. 

We write its simple normal crossing boundary $\cD$ as the union of (the proper transform of) the toric part $\cD_{\T}$ and the non-toric part $\cD\setminus \cD_{\T}$. Note that, by property $(T')$, there is a canonical embedding of the set of torus fixed points $\cD^{[0]}_{\T}$ of $\barcT_{\T}$ into the zero-dimensional stratum $\cD^{[0]}$.  

This completes the proof of the Proposition.
\end{proof}
 
\begin{definition}\label{kExtensions} Let $\cW^{-1}_{n, \eps = k^{-1}}\overline{\tilde{\alpha}_{\cX}(k\cL)}, \cW^{-1}\tilde{\beta}_{\cX}(k\cL)$ denote the classes of holomorphic $n+1$-forms corresponding to the multiple $k \cL$ for $k \gg 1$, \emph{where the deformation parameter $\eps$ is specialised to $\eps = k^{-1}$}. We denote by $[\tilde{\omega}^{(k)}_{\pm}]$ the classes of logarithmic extensions of $\cW^{-1}_{n, \eps = k^{-1}}\overline{\tilde{\alpha}_{\cX}(k\cL)}, \cW^{-1}\tilde{\beta}_{\cX}(k\cL)$ to $\barcT$ provided by Proposition \ref{GoodCompactProp}, and by $[\omega^{(k)}_{\pm}]$ their specialisation to $z = 0$. 
\end{definition}
\begin{cor}\label{FutakiToCohoIFCor} Suppose that 
\begin{itemize} 
\item[$\dagger$] for $k \gg 1$, the compactifications $(\barcT_k, \nabla_{\cW_k})$ corresponding to $\cW_k := \cW(k\cL)$ constructed in Proposition \ref{GoodCompactProp} satisfy the stationary phase formula, i.e. condition $(iv)$ of Definition \ref{CompactPropertiesDef}. Note that, according to Lemma \ref{StatPhasLem}, a sufficient condition is that $\Crit(\cW_k)$ is discrete and we have
\begin{equation}\label{MorseCondition}
\sum_{p \in \Crit(\cW_k)} m_p = (-1)^{n+1} \chi(\barcT_k \setminus \cD_k).
\end{equation}
\end{itemize}
Then the logarithmic extensions $[\tilde{\omega}^{(k)}_{\pm}]$, $[\omega^{(k)}_{\pm}]$ of Definition \ref{kExtensions} satisfy
\begin{align*}
&\DF(\cX, \cL) = k^{-n}\DF(\cX, k\cL)= \left(\frac{1}{2\pi \ii z} \right)^{n+1}\bra k^{-n}\tilde{\omega}^{(k)}_-, \tilde{\omega}^{(k)}_+ \ket_{ch, \nabla^z_{\cW}} + O(k^{-1})\\
&  = \left( \frac{1}{2\pi \ii}\right)^{n+1} \bra k^{-n}\omega^{(k)}_-,\omega^{(k)}_+ \ket_{ch} + O(k^{-1}).
\end{align*}  
As a consequence, the residue formula
\begin{align}\label{DFResThm}
\nonumber&\DF(\cX, \cL) \\
&= (-1)^{n+1} \sum_{P_{n+1} : \cD_k(P_{n+1}) \neq \emptyset} \langle \res_{P_{n+1}}(\omega^{(k)}_+)\big| \res_{P_{n+1}}(\nabla_{\cW_k})^{-1}\big| \res_{P_{n+1}}(k^{-n}\omega^{(k)}_-)\rangle + O(k^{-1}),
\end{align}
holds, where the right hand side is computed on $\barcT_k$, as a sum over the zero-dimensional strata of the boundary $\cD_k$, and $\nabla_{\cW_k}$, $\omega^{(k)}_{\mp}$ are determined by the K\"ahler parameter $k\cL$. 

\end{cor} 
\begin{proof} The first identity $\DF(\cX, \cL) = k^{-n}\DF(\cX, k\cL)$ follows at once from \eqref{OdakaWangFormula}. The identities involving the cohomology intersection form follow from Lemma \ref{FutakiToCohoIF} after replacing $\cL$ with $k\cL$ and setting $\eps = k^{-1}$. Finally, the residue formula follows at once from Theorem \ref{MainResThm}.
\end{proof}
\begin{exm} In the case of the degeneration to the normal cone $(\cX, k\cL_r)$ in $\PP^1$ considered above, $V(\cW(k\cL_r)) \subset \barcT_{\T}$ is a smooth curve of genus $1$ and so $\{\cW(k\cL_r) = 0\} \subset \T$ is given by the complement of 5 points in a a genus $1$ curve (the intersections $V(\cW(k\cL_r)) \cap \cD_{\T}$), thus we have $\chi(\{\cW(k\cL_r) = 0\}) = -5$ and $\chi(\barcT_{\T} \setminus \cD) = \chi(\T \setminus \{\cW(k\cL_r) = 0\}) = 5$. On the other hand one can check directly that, for $r$ sufficiently close to $\frac{1}{2}$ (i.e. for $\cL_r$ sufficiently close to $-\frac{1}{2}K_{\cX}$) and for all $k \gg 1$, the LG potential $\cW(k\cL_r)$ has $5$ isolated critical points, not contained in the locus $\{\cW(k\cL_r) = 0\}$ (see e.g. \cite{Stoppa_LargeComplex}, Section 5.2). So in this example the condition $(M)$ is satisfied and the stationary phase formula holds.  
\end{exm}
\begin{exm} When $(\cX, \cL_r)$ is the product test configuration given by $(\PP(\olo\oplus\olo(1)) \cong \Bl_p \PP^2, H - r E)$, the locus $V(\cW(k\cL_r)) \subset \barcT_{\T} = \Bl_{p_3, p_5} S_6$ is a smooth curve of genus $1$ and so $\{\cW(k\cL_r) = 0\} \subset \T$ is given by the complement of 4 points in a a genus $1$ curve (the intersections $V(\cW(k\cL_r)) \cap \cD_{\T}$). So in this case we have $\chi(\{\cW(k\cL_r) = 0\}) = -4$, $\chi(\barcT_{\T} \setminus \cD) = \chi(\T \setminus \{\cW(k\cL_r) = 0\}) = 4$. One can also compute directly that for all $r \in (0, 1)$ and for all $k \gg 1$, the LG potential $\cW(k\cL_r)$ has $4$ isolated critical points, not contained in the locus $\{\cW(k\cL_r) = 0\}$ (see e.g. \cite{Stoppa_LargeComplex}, Example 3.1). Thus the condition $(M)$ is satisfied and the stationary phase formula holds. 
\end{exm}
\begin{exm} It is natural to ask if condition $(M)$ also holds for $\cW_{\cX}(-k K_{\cX})$, for $k \gg 1$, when $\cX = \Bl_{E \times \{0\}}(\Bl_p \PP^2 \times \PP^1)$, with respect to the compactification for $\T$ obtained from $\barcT_{\T} = \widetilde{T}_{(\cP')^{\circ}}$. However in this case a direct computation seems out of reach.
\end{exm}
\subsection{Residue theorem and localisation formulae}\label{MainProofSubSec}
Our aim in the present Section is to show that that the right hand side of the residue formula \eqref{DFResThm} can be interpreted naturally in terms of suitable Donaldson-Futaki invariants for a compactification $X^{\vee}$ of the torus $\opT \cong (\C^*)^n$. This requires several steps. 
\begin{prop}\label{FibrationsProp} There exists a finite set of toric fibrations $\pi_i\!:\cY^{\vee, i} \to \PP^1$, for $i = 1, \ldots, m$, such that
\begin{enumerate}
\item[$(i)$] $\cY^{\vee, i}$ is smooth, and the central fibre $\cY^{\vee, i}_0$ is a union of smooth irreducible divisors;
\item[$(ii)$] the general fibre of each $\pi_i$ is a fixed compactification $X^{\vee}$ of the structure torus $\opT$;

\item[$(iii)$] $\cY^{\vee, i}$ can be endowed with a (non-unique) toric vector field $v_i$, with isolated fixed points, which is compatible with $\pi_i$ in the sense that, after rescaling, $v_i$ generates a $\C^*$-action covering the standard $\C^*$-action on $\C \subset \PP^1$.
\item[$(iv)$] Let $Z(v_i)$ denote the set of fixed points of $v_i$ contained in the central fibre $\cY^{\vee, i}_0$, and let $Z_D(v_i) \subset Z(v_i)$ be the set of such fixed points which are contained in an irreducible divisor $D \subset \cX^{\vee, i}_0$. Then we can choose $v_i$ as above so that for all $D \subset \cX^{\vee, i}_0$ we have  
\begin{align}\label{RankInequality}
2 |Z(v_i)| \leq h^{1,1}(\cY^{\vee, i}) + \dim D^{\perp} +  n + 1+ Z_D(v_i),  
\end{align}   
where $D^{\perp} \subset H^{1,1}(\cY^{\vee, i})$ denotes the space of $(1, 1)$-classes restricting to zero on $D$.
\item[$(v)$] It is possible to choose a toric compactification $\barcT_{\T}$ of $\T$ as in Lemma \ref{IntermediateCompactLem}, with toric boundary $\cD_{\T}$ and torus fixed points $\cD^{[0]}_{\T} \subset \cD_{\T}$ (i.e. the zero-dimensional stratum), so that each set of fixed points $Z(v_i) \subset \cY^{\vee, i}_0$ contains a distinguished subset $\widetilde{Z}(v_i) \subset Z(v_i)$, which has a natural inclusion $\widetilde{Z}(v_i) \subset \cD^{[0]}_{\T}$, inducing a partition 
\begin{equation}\label{ZeroesDecomposition}
\cD^{[0]}_{\T} = \cup^m_{i=1} \widetilde{Z}(v_i).
\end{equation} 
\end{enumerate}
\end{prop}
\begin{rmk} The condition \eqref{RankInequality} will be used for the construction of the complex $(1,1)$-classes $[\eta^i_k]$, $[\xi^i_k]$ appearing in Theorem \ref{MainThmIntro}. The point is having some condition that controls (twice) the number of fixed points in terms of the K\"ahler rank, the dimension and the boundary, in order to show the surjectivity of an evaluation map for the two complex Hamiltonians corresponding to $[\eta^i_k]$ and twists $[\xi^i_k]$, which can be satisfied after ``breaking up" $\barcT_{\T}$ in sufficiently many, sufficiently simple pieces.  
\end{rmk}
\begin{proof} Let $\barcT_{\T}$ be a toric compactification of $\T$ constructed in Lemma \ref{IntermediateCompactLem}, with toric boundary $\cD_{\T}$ and torus fixed points $\cD^{[0]}_{\T} \subset \cD_{\T}$. Recall that $\barcT_{\T}$ was constructed from the datum of a a toric mirror LG potential  $\cW$ mirror to a regular toric test configuration $\cX$ for $X$. As we explained in Section \ref{DFResSec}, we may assume without loss of generality that there is a toric morphism $\cX \to X \times \PP^1$, corresponding to a decomposition $N_{\cX} \cong N \times \Z$. Thus, from the proof of Lemma \ref{IntermediateCompactLem}, we see that there is an induced toric morphism $\barcT_{\T} \to \PP^1$, and we denote a fixed smooth fibre by $X^{\vee}$. 
\begin{rmk}
Note that when $X = T_{P}$ is a toric Fano with reflexive polytope $P$, by the proof of Lemma \ref{IntermediateCompactLem}, we can choose $\barcT_{\T}$ so that $X^{\vee}$ is the polar dual $T_{P^{\circ}}$.
\end{rmk}   

Recall $\barcT_{\T}$ is not unique. Given a choice of $\barcT_{\T}$, we may replace it with a new smooth compactification, satisfying the properties of Lemma \ref{IntermediateCompactLem}, obtained by a refinement of the toric fan of the original $\barcT_{\T}$. Choosing the refinement to be trivial on $N_{\R} \times \{0\} \subset N_{\cX}$ does not change the general fibre $X^{\vee}$. Thus we can fix a choice of $\barcT_{\T}$ with the property that its fan can be decomposed into subfans, such that each fan in the decomposition corresponds to a smooth quasi-projective toric variety $\barcT^{i, o}_{\T}$ mapping to $\C$, for $i = 1, \ldots, m$, which can be compactified canonically in the standard way, by gluing with the trivial family at infinity, to a projective toric variety $\barcT^{i}_{\T}$ with a map $\barcT^{i}_{\T} \to \PP^1$, with general fibre $X^{\vee}$ and central fibre given by the union of smooth irreducible divisors. Moreover, we may choose this decomposition so that the sets of maximal cones are disjoint.

By refining our decomposition of the toric fan of $\barcT_{\T}$ (so, in particular, increasing $m$), we can ensure that the inequality \eqref{RankInequality} holds for all irreducible divisors $D \subset \cY^{\vee, i}_0$.

A generic rational choice of a toric vector field $v_i$ on $\barcT^{i}_{\T}$ is compatible with the map $\barcT^{i}_{\T} \to \PP^1$, in the sense that after rescaling it generates a $\C^*$-action covering the standard $\C^*$-action on $\C \subset \PP^1$, and has isolated fixed points on $\barcT^{i}_{\T}$. Moreover, by construction, the set of fixed points $Z(v_i)$ of $v_i$ in the central fibre $\barcT^{i}_{\T, 0}$ contains a distinguished subset $\widetilde{Z}(v_i) \subset Z(v_i)$ which corresponds to a set of torus fixed points on $\barcT_{\T}$, i.e. we have $\widetilde{Z}(v_i) \subset \cD^{[0]}_{\T}$ with $\widetilde{Z}(v_i) \cap \widetilde{Z}(v_j) = \emptyset$ for $i \neq j$, yielding a partition \eqref{ZeroesDecomposition}.  
\end{proof}
The following special case follows immediately.
\begin{cor}\label{SimpleFibrationsCor} In the setup of Proposition \ref{FibrationsProp}, if we do not require the condition \eqref{RankInequality}, then we can choose $m = 2$, and we denote the fibrations obtained as the canonical compactifications of $\barcT_{\T}|_{\PP^1\setminus \{\infty\}} \to \C$, $\barcT_{\T}|_{\PP^1\setminus \{0\}} \to \C$ (in the coordinates Proposition \ref{FibrationsProp}) by $\pi^0\!: \cY^{\vee, 0} \to \PP^1$, respectively $\pi^{\infty}\!:\cY^{\vee, \infty} \to \PP^1$. 
\end{cor}
We illustrate the construction of the fibrations $\pi_i\!: \cY^{\vee, i} \to \PP^1$ in some examples.
\begin{exm}\label{NormalConeAnticanTest} Recall that, for the degeneration to the normal cone of a point in $X = \PP^1$, we have $\barcT_{\T} = \Bl_{p_4} S_6$. Up to the action of $GL(2, \Z)$, the fan of $\barcT_{\T}$ is spanned by the vectors
\begin{align*}
& w_1 = (1,0),\,w_2 = (1,1),\,w_3 = (0,1),\\
& w_4 = (-1,1),\,w_5 = (-1,0),\,w_6 = (-1,-1),\,w_7=(0,-1),
\end{align*} 
with corresponding maximal cones $K_i$ spanned by $w_i, w_{i +1}$ for a cyclic index $i$. An admissible decomposition of the set of cones is given by $\{K_1, K_2, K_3, K_4\} \cup \{K_5, K_6, K_7\}$. The fibration $\cY^{\vee, 1}$ for $X^{\vee} \cong \PP^1$ corresponding to $\{K_1, K_2, K_3, K_4\}$ is isomorphic to $\Bl_{q_1 \times \{0\}, q_2 \times \{0\}}\PP^1\times\PP^1$. It has central fibre $\cY^{\vee, 1}_0$ given by a tree of $3$ smooth rational curves, containing $4$ torus fixed points. Let $D_1$ be any irreducible component of $\cY^{\vee, 1}_0$. Let $v_1$ be a generic toric vector field on $\cY^{\vee, 1}$. Then we have $ h^{1,1}(\cY^{\vee, 1}) + \dim D^{\perp}_1 + n +1 \geq 7 > 6 = 8 - 2 = 2|Z(v_1)| - |Z_{D_1}(v_1)|$. Similarly, $\cY^{\vee, 2}$, corresponding to $\{K_5, K_6, K_7\}$, is isomorphic to $\Bl_{q \times \{0\}} \PP^1 \times \PP^1$ and has central fibre given by the transverse intersection of $2$ smooth rational curves, containing $3$ torus fixed points, so for a generic toric vector field $v_2$ on $\cY^{\vee, 1}$ we have $h^{1,1}(\cY^{\vee, 2}) + \dim D^{\perp}_2+ n+1 \geq 6 > 6 - 2 = 2|Z(v_2)| -|Z_{D_2}(v_2)|$ for any irreducible component $D_2$ of $\cY^{\vee, 2}_0$.
\end{exm}
\begin{exm}\label{ProductAnticanTest} Recall that, for the product test configuration $\cX \cong \PP(\olo\oplus\olo(1))$, we have $\barcT_{\T} = \Bl_{p_3, p_5} S_6$. Up to the action of $GL(2, \Z)$, the fan of $\barcT_{\T}$ is spanned by the vectors
\begin{align*}
& w_1 = (1,0),\,w_2 = (1,1),\,w_3 = (1,2),\,w_4 = (0,1),\\
& w_5 = (-1,0),\,w_6 = (-2, -1),\,w_7 = (-1,-1),\,w_8=(0,-1),
\end{align*} 
with corresponding maximal cones $K_i$ spanned by $w_i, w_{i +1}$ for a cyclic index $i$. An admissible decomposition of the set of cones is given by $\{K_1, K_2, K_3, K_4\} \cup \{K_5, K_6, K_7, K_8\}$. The corresponding $\cY^{\vee, i}$ for $X^{\vee} \cong \PP^1$, for $i =1, 2$, are isomorphic to $\Bl_{q' \times \{0\}, q'' \times \{0\}}\PP^1 \times \PP^1$ where $q' \times \{0\}$, $q'' \times \{0\}$ are infinitely close. The central fibre $\cY^{\vee, i}_0$ is given by a tree of $3$ smooth rational curves, containing $4$ torus fixed points, so for a generic toric vector field $v_i$ we have $h^{1,1}(\cY^{\vee, i}) + \dim D^{\perp}_i + n + 1 \geq 8 > 8 - 2 = 2|Z(v_i)| - |Z_{D_i}(v_i)|$ for any irreducible component $D_i \subset \cY^{\vee, i}_0$.  
\end{exm}
\begin{exm}\label{SlopeUnstableSurfaceExm} We showed that, when $(\cX, -K_{\cX})$ is the Fano toric test configuration for $X = \Bl_p \PP^2$ given by the degeneration to the normal cone of $E$, $\cX = \Bl_{E \times \{0\}} X \times \PP^1$, then an admissible choice of intermediate toric compactification for $\T$ is given by $\barcT_{\T} = \widetilde{T}_{(\cP')^{\circ}}$, a crepant resolution of the Gorenstein canonical Fano threefold $T_{(\cP')^{\circ}}$ corresponding to some maximal triangulation. 

Note that $(\cP')^{\circ} \cap \{z = 0 \}$ is the reflexive polygon
\begin{equation*}
P^{\circ} = \conv(\ell),\,\ell^{\circ}:=\{( 2, -1), (-1,  2), (-1,  0), ( 0, -1)\}
\end{equation*}
polar dual to
\begin{equation*}
P = \conv(\{( 1,  0), ( 0,  1), ( 1,  1), (-1, -1) \}),
\end{equation*}
and so in our construction above each $\cY^{\vee, i}$ gives a toric fibration with generic fibre $X^{\vee} = \widetilde{T}_{P^{\circ}}$, the crepant resolution of the orbifold $T_{P^{\circ}}$ given by a maximal triangulation of $P^{\circ}$, which is (after resolving the base locus of $W_X(-K_X)$) the mirror of $X = T_{P}$. One can check that $T_{P^{\circ}}$ is the singular intersection of two quadrics in $\PP^4$, 
\begin{equation*}
T_{P^{\circ}} = \{x_1 x_3 = x_0 x_4,\, x_2 x_4 = x^2_0\} \subset \PP[x_0 : x_1 : x_2 : x_3 : x_4],
\end{equation*}
and its resolution is $X^{\vee} \cong \Bl_{q'_2, q'_3} S_6 = \Bl_{q'_2, q'_3} \Bl_{q_1, q_2, q_3} \PP^2$, where $q'_2$, $q'_3$ are points lying on the exceptional divisors over $q_2$, $q_3$. 

On the other hand, because of the large number of cones appearing in the fan of $\barcT_{\T} = \widetilde{T}_{(\cP')^{\circ}}$, it is difficult to describe explicitly a collection $\cY^{\vee, i}$, for $i =1, \ldots, m$, such that \eqref{RankInequality} holds for all $\cY^{\vee, i}$. 

Let as assume for a moment that we can work with the Gorenstein canonical compactification $T_{(\cP')^{\circ}}$ rather than its crepant resolution. 

We note that there is a presentation
\begin{equation*}
(\cP')^{\circ} = \conv(\tilde{\ell}^{\circ} \times \{z = -1\} \cup \ell^{\circ}\times \{z = 0\} \cup \ell^{\circ}\times \{z = 1\}),
\end{equation*}
where $\tilde{\ell}^{\circ} := \{( 2, -1),(-1,  2),( 1, -1),(-1,  1)\}$. Then, at least formally, we can choose our fibrations \emph{for the orbifold} $X^{\vee}$ as $\cY^{\vee, i} = T_{\cQ^i}$, where
\begin{align*}
& \cQ^{1} := \conv((0,0,-1) \cup \ell^{\circ}\times \{z=0\} \cup \ell^{\circ}\times \{z=1\} \setminus \{( 2, -1,1),( 0, -1,1)\}),\\
& \cQ^{2} := \conv((0,0,-1) \cup \ell^{\circ}\times \{z=0\} \cup \ell^{\circ}\times \{z=1\} \setminus \{(-1,  2,1),(-1,  0,1)\}),\\
& \cQ^{3} := \conv((0, 0, 1) \cup \tilde{\ell}^{\circ} \times \{z = -1\} \cup \ell^{\circ}\times \{z=0\} \setminus \{( 2, -1, -1),( 1, -1, -1)\}),\\
&\cQ^{4} := \conv((0, 0, 1) \cup \tilde{\ell}^{\circ} \times \{z = -1\} \cup \ell^{\circ}\times \{z=0\} \setminus \{(-1,  2, -1),(-1,  1, -1)\}).
\end{align*}  
Indeed, in each case, we can compute that the Picard rank of $\cY^{\vee, i}$ is $3$ and the number of torus fixed points on the central fibre is $5$, from which \eqref{RankInequality} follows. Note that there is a reflection symmetry along the plane $x = y$ yielding isomorphisms $\cY^{\vee, 1} \cong \cY^{\vee, 2}$, $\cY^{\vee, 3} \cong \cY^{\vee, 4}$. 
\end{exm}

\begin{lemma}\label{DFDecompositionLem} Let $\cW_k := \cW(k\cL)$, with corresponding linear mirror isomorphism 
\begin{equation*}
\Theta_{k\cL}|_{|z=0}\!: \Jac(\cW_k) \xrightarrow{\,\,\cong\,\,} (H^*(\cX, \C), *_{[c_1(k\cL)]})
\end{equation*}
(see property $(ii)$ of Section \ref{ToricLGSec}). Define
\begin{equation}\label{ThetaFunctionsLargeK}
[\theta_k] = k^{-1}\Theta^{-1}_{k\cL}(c_1(k\cL))|_{z=0},\, [\psi_k] = \Theta^{-1}_{k\cL}(c_1(K_{\cX/\PP^1}))|_{z=0} + [\cW_k] \in  \Jac(\cW(k\cL)).
\end{equation} 
Suppose that the stationary phase formula (condition $\dagger$ in Corollary \ref{FutakiToCohoIFCor}) holds for $(\barcT_k, \nabla_{\cW_k})$. Then, for any choice of fibrations $\cY^{\vee, i}$, $i =1, \ldots, m$ as in Proposition \ref{FibrationsProp}, we have, in terms of the classes $[\omega^{(k)}_{\pm}]$ of Definition \ref{kExtensions},
\begin{align}\label{DFDecomposition}  
\nonumber &\DF(\cX, \cL) 
 = (-1)^{n+1} \sum^m_{i=1}\sum_{ \cD(P_{n+1}) = p \in \widetilde{Z}(v_i) }  [\cW^{-1}_k\theta_k]^n\left(\frac{n c}{n+1}[\cW^{-1}_k\theta_k] - 1 + [\cW^{-1}_k\psi_k] \right)\big|_{p} \\
\nonumber &\quad\quad\quad\quad\quad\bra\res(\nabla_{\cW(k\cL)})\big|_{p}  \Omega_0|_p, \Omega_0|_p\ket\left(1 + O(k^{-1})\right)\\
& + (-1)^{n+1} \sum_{\cD(P_{n+1}) \in \cD \setminus \cD^{[0]}_{\T}}   \langle \res_{P_{n+1}}(\omega^{(k)}_+)\big| \res_{P_{n+1}}(\nabla_{\cW(k\cL)})^{-1}\big| \res_{P_{n+1}}(k^{-n}\omega^{(k)}_-)\rangle + O(k^{-1}),
\end{align}
where evaluation for classes in $\opH^n_{\mp}$ is defined in terms of restrictions along toric strata (i.e. in terms of intersections $\cD(P_{n+1})$), as explained in the discussion of the global residue theorem, Section  \ref{ResThmSubSec}. Note that the subsets $\widetilde{Z}(v_i) \subset Z(v_i) \subset \cY^{\vee, i}_0$ appear in the decomposition.
\end{lemma}
\begin{proof} Recall that the classes of logarithmic forms $ [\omega^{(k)}_{\pm}]$ are given explicitly by 
\begin{align*}
&[\omega^{(k)}_-] =  \cW^{-1}_{n, k^{-1}}(k\cL) \Theta^{-1}_{k\cL} \left((c_1(k\cL))^n\right)\big|_{z=0},\\
&[\omega^{(k)}_+] = \cW^{-1}(k\cL) \Theta^{-1}_{k\cL} \left(\frac{n (k^{-1}c)}{n+1}c_1(k\cL) + c_1(K_{\cX / \PP^1})\right)\big|_{z=0}. 
\end{align*}
As in the proof of Proposition \ref{GoodCompactProp} we note that
\begin{align*}
& k^{-n}\Theta^{-1}_{k\cL} \left((c_1(k\cL))^n\right)\big|_{z=0} = k^{-n}\Theta^{-1}_{k\cL}\big(c_1(k\cL) *_{[c_1(k\cL)]} \cdots *_{[c_1(k\cL)]}  c_1(k\cL) +O(k^{-1}) )\\
& = \left(k^{-1}\Theta^{-1}_{k\cL}(c_1(k\cL))\right)^n  + O(k^{-1}).
\end{align*}
The result now follows at once from \eqref{DFResThm}, by labelling the contributions to the residue formula according to \eqref{ZeroesDecomposition}.
\end{proof}

\subsection{Completion of the proofs of Theorems \ref{MainThmIntro}, \ref{VanishingThm} and \ref{MainThmVariant}} \label{MainProofSubSecCompletion}
In the present Section we will complete the proofs of Theorems \ref{MainThmIntro}, \ref{VanishingThm} and \ref{MainThmVariant}.

\begin{definition}\label{TestConfigDef} For each $i = 1, \ldots, m$, we let $(\cX^{\vee, i}, v_i, \eta^i, \xi^i)$ denote a regular toric compactified test configuration for $(X^{\vee}, \eta^i|_{X^{\vee}}, \xi^i|_{X^{\vee}})$ obtained by endowing the total space of the fibration $\pi^i\!:\cY^{\vee, i} \to \PP^1$ with any vector field $v_i$ as in Proposition \ref{FibrationsProp}, $(iii)$ and with $(1,1)$-classes $[ \eta^i]$, $[\xi^i]$. This is \emph{formal} in the sense that $[ \eta^i]$, $[\xi^i]$ are complex and are not required to satisfy positivity conditions. 

When the construction of Proposition \ref{FibrationsProp} is performed with respect to the LG potential $\cW(k\cL)$ for $k \gg 1$, we allow the classes $(\eta^i, \xi^i)$ to depend on $k$ and sometimes write $(\eta^i_k, \xi^i_k)$ to emphasise this dependence.  

Similarly, in the situation of Corollary \ref{SimpleFibrationsCor}, when we do not require the condition \eqref{RankInequality}, then we denote by $(\cX^{\vee, 0}, v_0, \cZ^0_k, \tilde{\cZ}^{0}_k)$,  $(\cX^{\vee, \infty}, v_{\infty}, \cZ^{\infty}_k, \tilde{\cZ}^{\infty}_k)$ the regular toric compactified test configurations for $(X^{\vee}, \cZ^0_k|_{\cX^{\vee, 0}_1}, \tilde{\cZ}^{0}_k|_{\cX^{\vee, 0}_1})$, $(X^{\vee}, \cZ^0_k|_{\cX^{\vee, 0}_1}, \tilde{\cZ}^{\infty}_k|_{\cX^{\vee, \infty}_1})$, in the sense of Definition \ref{DFExtendedKahler},  obtained by endowing the total space of the fibrations $\pi^0\!:\cY^{\vee, 0} \to \PP^1$, $\pi^{\infty}\!:\cY^{\vee, 0} \to \PP^1$ with any vector field $v_i$ as in Proposition \ref{FibrationsProp}, $(iii)$ and with extended K\"ahler parameters $\cZ^0_k, \tilde{\cZ}^0_k \in H^*(\cX^{\vee, 0}, \C)$, $\cZ^{\infty}_k, \tilde{\cZ}^{\infty}_k \in H^*(\cX^{\vee, \infty}, \C)$.     
\end{definition}
\begin{proof}[Proof of Theorem \ref{MainThmIntro}] Suppose $(\cX^{\vee}, v, \eta, \xi):=(\cX^{\vee, i}, v_i, \eta^i_k, \xi^i_k)$ is one of the test configurations of Definition \ref{TestConfigDef}. We claim that it is possible to choose the  holomorphic vector field $v_i$ (as in $(iii)$ of Proposition \ref{FibrationsProp}) and the $(1,1)$-classes $[\eta^i_k]$, $[\xi^i_k]$ such that, in the large volume limit (i.e. up to $O(k^{-1})$), the contribution 
 \begin{align}\label{LocalContrib}
(-1)^{n+1} \sum_{ \cD(P_{n+1}) = p \in \widetilde{Z}(v_i) }  [\cW^{-1}_k\theta_k]^n\left(\frac{n c}{n+1}[\cW^{-1}_k\theta_k] - 1 + [\cW^{-1}_k\psi_k] \right)\bra\res(\nabla_{\cW(k\cL)})   \Omega_0 , \Omega_0\ket \big|_{p}
\end{align}
appearing in \eqref{DFDecomposition} can be identified naturally with the Atiyah-Bott equivariant localisation formula \eqref{DFAtiyahBottLocalisation} for $\DF(\cX^{\vee, i}, v_i, [\eta^i_k], [\xi^i_k])$, namely
\begin{align*} 
&\nonumber \DF(\cX^{\vee, i}, [\eta^i_k], [\xi^i_k])
= \sum_{p \in Z(v_i)} \frac{(-h_{\eta}(p))^{n}\big(-\frac{nc^{\vee}_{\eta, \xi}}{n+1} h_{\eta}(p)-\sum^{n+1}_{i=0} w_i(p) +1 -  h_{\xi}(p)\big)}{e(T_p)(v)}.  
\end{align*}
We emphasise that determining the (generic) vector field $v_i$ is also part of the problem. 

Theorem \ref{MainThmIntro} follows immediately from this claim, by defining the base loci contributions as
\begin{equation}\label{BaseLocusTerm}
\cB(\cW_k) := (-1)^{n+1} \sum_{\cD(P_{n+1}) \in \cD \setminus \cD^{[0]}_{\T}}   \langle \res_{P_{n+1}}(\omega^{(k)}_+)\big| \res_{P_{n+1}}(\nabla_{\cW(k\cL)})^{-1}\big| \res_{P_{n+1}}(\omega^{(k)}_-)\rangle.
\end{equation}   
\begin{rmk}\label{NoncompactDFRmk} Applying the same large volume limit expansion appearing in the proof of Lemma \ref{DFDecompositionLem} to the Grothendiek residue expression \eqref{DFGrothRes}, when the critical points of $\cW(k\cL)$ are nondegenerate, we find immediately
\begin{equation*}
\DF(\cX, \cL) = \sum_{p \in \Crit(\cW(k\cL))} \frac{(\theta_k)^n(\frac{n c}{n+1}\theta_k - \cW_k + \psi_k)}{\prod^{n+1}_{i=1} x^2_i\det \nabla^2\cW_k}\big|_p + O(k^{-1}).
\end{equation*} 
The leading term clearly resembles an \emph{ill-defined} localised Donaldson-Futaki invariant on the \emph{noncompact} manifold $\T$, with respect to the vector field $\nabla \cW_k$. Our construction aims precisely at turning this into a sum of well-defined formal Donadson-Futaki invariants for compactified mirrors. 
\end{rmk}

Let us write
\begin{align*}
&\tilde{\theta} := \cW^{-1}_k \theta_k,\,\tilde{\psi} := \cW^{-1}_k\psi_k,\\
& f^{n+1}(p) := |P_{n+1}\!:\cD(P_{n+1}) = p \in \cD_{\T}|\bra\res(\nabla_{\cW(k\cL)})\big|_{p}  \Omega_0|_p, \Omega_0|_p\ket 
\end{align*}  
where $\theta_k$, $\psi_k$ are defined in \eqref{ThetaFunctionsLargeK} (note that we omit $k$ in the notation for simplicity, but all the quantities we consider in the following depend on $k$).  Then the contribution \eqref{LocalContrib} can be written as 
\begin{align*}
\sum_{ p \in \widetilde{Z}(v_i)} \left(-  f\tilde{\theta} \right)^n \left(\frac{n c}{n+1}\left( - f \tilde{\theta} \right) + f \left(1 - \tilde{\psi} \right)  \right) \big|_{p}. 
\end{align*}
If $p$ is a zero of a holomorphic vector field $v$, we set
\begin{equation*}
d(v)|_p = (\det(\nabla v))^{\frac{1}{n+1}}|_p,\,t(v)|_p = \operatorname{tr}(\nabla v)|_p-1. 
\end{equation*}
These quantities are well defined for all $v \in \mathfrak{t} := \operatorname{Lie}(\T)$, not necessarily generating a $\C^*$-action, and agree with the corresponding expressions involving the equivariant Euler class and sum of weights in the case of $\C^*$-actions.

Let us introduce the functions $H_{i}\!: Z(v_i) \to \C$, $K_{i}\!: Z(v_i) \to \C$, such that
\begin{align*}
& H_{i}|_{\widetilde{Z}(v_i)} := d(v_i) f\tilde{\theta},\,\,H_{i}|_{Z(v_i)\setminus\widetilde{Z}(v_i)}:= 0,\\
& K_{i}|_{\widetilde{Z}(v_i)} := - t(v_i) - d(v_i) f \left(1 - \tilde{\psi} \right),\,\,K_{i}|_{Z(v_i)\setminus\widetilde{Z}(v_i)}:= 0.   
\end{align*}

Then, we may write the contribution \eqref{LocalContrib} in the form
\begin{align*}
&\sum_{ p \in \widetilde{Z}(v_i)} \frac{\left(- d(v_i) f\tilde{\theta} \right)^n \left(\frac{n c}{n+1}\left( -  d(v_i) f \tilde{\theta} \right) -t(v_i)+ \left(t(v_i) + d(v_i) f \left(1 - \tilde{\psi} \right)  \right)\right)}{(d(v_i))^{n+1}} \big|_{p}\\
& = \sum_{ p \in Z(v_i)} \frac{\left(- H_i \right)^n \left(-\frac{n c}{n+1} H_i  - t(v_i) -  K_i \right)}{(d(v_i))^{n+1}} \big|_{p}. 
\end{align*}

Thus our claim follows if we can realise $H_i$, $K_i$ as the values of \emph{complexified} Hamiltonians $h_{\eta^i}$, $h_{\xi^i}$ for $v_i$, i.e. determining \emph{both} a vector field $v_i$ as above and complex $(1, 1)$-forms $\eta^i$, $\xi^i$ such that 
\begin{align}\label{PrescribedH}
&\nonumber\iota_{v_i} \eta^i = \delbar h_{\eta^i},\,h_{\eta^i}(p) = H_i(p) + O(k^{-1}), \\
&\iota_{v_i} \xi^i = \delbar h_{\xi^i},\,h_{\xi^i}(p) = K_i(p) + O(k^{-1}), \textrm{ for all } p \in Z(v_i).
\end{align}

Recall that the values $h_{\eta^i}(p)$,$h_{\xi^i}(p)$ only depend on the cohomology class $[\eta^i]$, $[\xi^i]$, up to an overall constant. Thus, we can consider the corresponding problem in cohomology, so the set of values of $h_{\eta^i}(p)$,$h_{\xi^i}(p)$ as $[\eta^i]$, $[\xi^i]$ vary in $H^{1,1}(\cX^{\vee, i}, \C)^{\oplus 2}$ is a complex linear subspace of $(\C^{|Z(v_i)|})^{\oplus 2}$ of dimension $2 h^{1,1}(\cX^{\vee, i})$. Note that, since by construction $\cX^{\vee, i}$ is smooth, toric and projective, the existence of such Hamiltonians is automatic.

We observe that, by the definition of $\tilde{\psi}$ and our construction of the test configurations $\cX^{\vee, i}$ through Proposition \ref{FibrationsProp}, $\tilde{\psi}$ induces a rational function on $\cX^{\vee, i}$ for which there exists an irreducible component $D_i \subset \cX^{\vee, i}_0$ such that
\begin{equation*}
\tilde{\psi}|_{D_i} = 1 + O(k^{-1}),
\end{equation*}
and so we have
\begin{equation*}
K_{i}|_{\widetilde{Z}(v_i) \cap D_i}  = -t(v_i) + O(k^{-1}) = -\operatorname{tr}(\nabla v_i) + 1 + O(k^{-1}).
\end{equation*}
For $p \in D_i$, the quantity $\operatorname{tr}(\nabla v_i)|_p$ equals the value at $p$ of a Hamiltonian on $D_i$ with respect to the class $-c_1(K_{D_i}) = -c_1(\olo(D_i)|_{D_i})$ (see e.g. \cite{Tian_libro}, Section 3.2), up to an overall constant. Thus, twists $[\xi^i]$ of the form $c_1(\olo(D_i)) + [\tilde{\xi}^i]$ for $[\tilde{\xi}^i] \in D^{\perp}_{i}$, achieve the correct values along $D_i$, and letting $[\eta^i]$, $[\tilde{\xi}^i]$ vary in $H^{1,1}(\cX^{\vee, i}, \C)$, respectively $D^{\perp}_i$, yields a subspace of 
\begin{equation*}
\cV_i : = \C^{|Z(v_i)|} \oplus \C^{|Z(v_i)| - |Z_{D_i}(v_i)|}
\end{equation*}
of dimension $h^{1,1}(\cX^{\vee, i}) + \dim D^{\perp}_i$ (recall that $|Z_{D_i}(v_i)|$ denotes the number of isolated fixed points contained in $D_i$, as in Proposition \ref{FibrationsProp}).

Now we also let the holomorphic vector field $v_i$ vary in the generic locus in $\mathfrak{t}$ of the maximal torus. As $v_i$, $[\eta^i]$, $[\tilde{\xi}^i]$ vary, the target complex vector space $\cV_i$ remains fixed, and the values of $v_i$, $h_{\eta^i}(p)$, $h_{\xi^i}(p)$ yield a complex subspace $\mathfrak{L} \subset \mathfrak{t} \times \cV_i$ of dimension $h^{1,1}(\cX^{\vee, i}) + \dim D^{\perp}_i +  n + 1$.  At the same time, the values of $H_i$, $K_i$ as $v_i$ varies in the generic locus give an analytic subvariety $\mathfrak{M} \subset \mathfrak{t} \times \cV_i$ of dimension $n + 1$. By \eqref{RankInequality}, possibly after a small linear perburbation of $\mathfrak{L}$ of order $O(k^{-1})$, there exist a vector field $v_i$ and classes $[\eta^i]$, $[\xi^i]$ lying in the intersection $\mathfrak{L} \cap \mathfrak{M} \neq \emptyset$ and so satisfying our condition \eqref{PrescribedH}. 

This particular $v_i$ might not generate a $\C^*$-action, but it can be approximated by $\C^*$-actions, and this approximation only introduces an error term of order $O(k^{-1})$. 

The requirement $\int_{X^{\vee}} (\eta^i)^n \neq 0$ might not be satisfied for this particular  class $[\eta^i]$, but if that happens we can perturb it to $[\eta^i + k^{-1} \varphi]$ where $[\varphi]$ is a $(1,1)$-class on $\cX^{\vee, i}$ such that $\int_{X^{\vee}}  (\eta^i)^{n-1} \wedge \varphi \neq 0$.  
 
Finally, the equality 
\begin{equation*}
c = c^{\vee}_{\eta^i, \xi^i}
\end{equation*} 
can be achieved by scaling $[\eta^i]$, $[\xi^i]$ suitably (so that \eqref{PrescribedH} remains solvable).
  
The upshot is that we have an identity
\begin{align*}
&\sum_{ p \in Z(v_i)} \frac{\left(- H_i \right)^n \left(-\frac{n c}{n+1} H_i  - t(v_i) -  K_i \right)}{(d(v_i))^{n+1}} \big|_{p}\\
& = \sum_{p \in Z(v_i)} \frac{(-h_{\eta^i}(p))^{n}\big(-\frac{n c^{\vee}_{\eta^i, \xi^i}}{n+1} h_{\eta^i}(p) - \sum^{n+1}_{j=0} w_j(p) +1 -  h_{\xi^i}(p)\big)}{e(T_p)(v)} + O(k^{-1}),   
\end{align*}
as required. This completes the proof of our claim and so of Theorem \ref{MainThmIntro}. 
\end{proof}
\begin{exm} Suppose $(\cX, \cL_r)$ is given by the degeneration to the normal cone of a point in $\PP^1$. We know that in this case $\barcT_{\T} = \Bl_{p_4} S_6$, and the torus fixed points on the central fibre of the test configurations $\cX^{\vee, 1} \cong \Bl_{q_1 \times \{0\}, q_2 \times \{0\}} \PP^1 \times \PP^1$, $\cX^{\vee, 2} \cong \Bl_{q\times\{0\}} \PP^1 \times \PP^1$ are given by $\{p_2, p_3, p'_4, p''_4\}$ (with $p'_4$, $p''_4$ mapping to $p_4$), respectively $\{p_1, p_5, p_6\}$. In Example \ref{NormalConeAnticanExm} below we will compute the quantities $\tilde{\theta}$, $\tilde{\psi}$, $f$ in the case when $(\cX, \cL_r) \cong (\cX, -\frac{1}{2} K_{\cX})$ has parameter $r = \frac{1}{2}$. Using that computation, we find that the required values at the torus fixed points are achieved if
\begin{align*}
& H_1(p) = \frac{1}{2} d(v_1)|_p,\,p \in \{p_2, p_3, p'_4, p''_4\},\,H_2(p) = \frac{1}{2} d(v_2)|_p,\,p \in \{p_1, p_5, p_6\},\\
& K_1(p) = t(v_1)\big|_p,\, p \in \{p_3, p'_{4}\},\,K_1(p) = (t(v_1) + d(v_1))\big|_p,\, p \in \{p_2, p''_{4}\},\\  
&K_2(p) = t(v_2)\big|_p,\,p \in \{p_1, p_6\},\,K_2(p) = (t(v_2) + d(v_2))\big|_p,\, p = p_5.
\end{align*}
Let us show how these prescribed values can be achieved for $\cX^{\vee, 1}$ (i.e. solving for $v_1$, $[\eta_{\cX^{\vee, 1}}]$ and $[\xi_{\cX^{\vee, 1}}]$). A completely analogous computation holds for $\cX^{\vee, 2}$. 

A vector field $v_1$ on $\cX^{\vee, 1}$ is induced by $\hat{v}_1$ on $\PP^1\times\PP^1$. In suitable affine coordinates $(z, w)$ on $\PP^1 \times \PP^1$, under the natural identifications, we have 
\begin{equation*}
\hat{v}_1 = a  z \frac{\del}{\del z} + b w \frac{\del}{\del w}, 
\end{equation*}
and so, passing to the blowup with local coordinates $(z, \xi)$ or $(\eta, w)$, we find
\begin{align*}
& v_1 = (a + b ) z \frac{\del}{\del z} + b \xi \frac{\del}{\del \xi},\textrm{ respectively } v_1 = a \eta \frac{\del}{\del \eta} + (a  + b ) w \frac{\del}{\del w}.
\end{align*} 
Using this, we can compute
\begin{align*}
& t(v_1)|_{p_2} = 2 a + b - 1,\, t(v_1)|_{p'_4} =  a + 2 b - 1,\\
& t(v_1)|_{p_3} = -2 a - b - 1,\,t(v_1)|_{p''_4} =  - a - 2 b - 1\\
& d(v_1)|_{p_2} = a (a + b ),\,d(v_1)|_{p'_4} = (a + b ) b ,\\ 
& d(v_1)|_{p_3} = a (a + b ),\,d(v_1)|_{p''_4} = (a + b ) b. 
\end{align*}
Choosing $a = b = \frac{1}{3}$, our prescribed values on $\cX^{\vee, 1}$ become
\begin{align*}
& H_1(p) = \frac{2}{9},\,p \in \{p_2, p_3, p'_4, p''_4\},\\ 
& K_1(p) = 0,\, p \in \{p_2, p'_{4}\},\,K_1(p) = -2,\, p \in \{p_3, p''_{4}\}. 
\end{align*}
These can be satisfied by choosing 
\begin{align*}
\eta_1 = \pi^* p^*_1 \omega^{(1)}_{\PP^1},\,\xi_1 = \pi^* p^*_2 \omega^{(2)}_{\PP^1}, 
\end{align*}
where $\omega^{(i)}_{\PP^1}$ are suitable multiples of the Fubini-Study form, $\pi$ denotes the blow-down map to $\PP^1\times\PP^1$ and $p_i\! : \PP^1\times\PP^1 \to \PP^1$ are the projections.   
\end{exm}
\begin{exm} Suppose $\cX \cong \PP(\olo \oplus \olo(1))$ is a product test configuration for $\PP^1$. We know that in this case $\barcT_{\T} = \Bl_{p_3, p_5} S_6$, and the torus fixed points on the central fibres of the test configurations $\cX^{\vee, i} \cong \Bl_{q' \times \{0\}, q'' \times \{0\}}\PP^1 \times \PP^1$, for $i = 1, 2$, are given by $\{p_2, p'_3, p''_3, p_4\}$ (with $p'_3$, $p''_3$ mapping to $p_3$), respectively $\{p_1, p'_5, p''_5, p_6\}$ (with $p'_5$, $p''_5$ mapping to $p_5$). In Example \ref{ProductAnticanExm} below we will compute the quantities $\tilde{\theta}$, $\tilde{\psi}$, $f$ in the case when $(\cX, \cL_r) \cong (\cX, -\frac{1}{3} K_{\cX})$ has parameter $r = \frac{1}{3}$. Using that computation, we find that the required values at the torus fixed points are achieved if
\begin{align*}
& H_1(p) = \frac{1}{3} d(v_1)|_p,\,p \in \{ p_2, p'_3, p''_3, p_4\},\, H_2(p) = \frac{1}{3} d(v_2)|_p,\,p \in \{p_1, p'_5, p''_5, p_6 \},\\
& K_1(p) = t(v_1)\big|_p,\, p \in \{p_2, p'_{3}\},\, K_1(p) = (t(v_1) + d(v_1))\big|_p,\, p \in \{p''_{3}, p_4\},\\  
& K_2(p) = t(v_2)\big|_p,\,p \in \{p'_5, p_6\},\, K_2(p) = (t(v_2) + d(v_2))\big|_p,\, p = \{p_1, p''_5\}.
\end{align*}
Let us work on $\cX^{\vee, 1}$. On the blowup, with local coordinates $(z, \xi)$ or $(z, \eta)$, we find
\begin{align*}
& v_1 = (a + b ) z \frac{\del}{\del z} + b \xi \frac{\del}{\del \xi},\textrm{ respectively } v_1 = (a + 2 b ) z \frac{\del}{\del z} + b \eta \frac{\del}{\del \eta},
\end{align*}
from which
\begin{align*}
& t(v_1)|_{p_2} = a + 2b - 1,\, t(v_1)|_{p'_3} =  a + 3 b - 1,\\
& t(v_1)|_{p_4} = -  a -  b - 1,\,t(v_1)|_{p''_3} =  - a - 3 b - 1\\
& d(v_1)|_{p_2} = (a + b ) b,\,d(v_1)|_{p'_3} = (a + 2 b ) b ,\\ 
& d(v_1)|_{p_4} =  a b,\,d(v_1)|_{p''_3} = (a + 2 b ) b. 
\end{align*}
Choosing $a = 1, b = 0$, our prescribed values on $\cX^{\vee, 1}$ become
\begin{align*}
& H_1(p) = 0,\,p \in \{p_2, p'_3, p''_3, p_4\},\\ 
& K_1(p) = 0,\, p \in \{p_2, p'_{3}\},\,K_1(p) = -2,\, p \in \{p''_3, p_{4}\}. 
\end{align*}
These can be satisfied with the same chooices 
$\eta_1 = \pi^* p^*_1 \omega^{(1)}_{\PP^1},\,\xi_1 = \pi^* p^*_2 \omega^{(2)}_{\PP^1}$  
as in our previous example.
\end{exm}

\begin{proof}[Proof of Theorem \ref{VanishingThm}] Suppose $\cA$ is an ample line bundle on $\cX$, such that
\begin{equation*}
\Theta_{k \cA }(\cA)|_{z=0} = r \cW(k \cA ) + O(k^{-1})
\end{equation*}
for some $r > 0$. This happens e.g. if $\cX$ is Fano and $\cA = - r K_{\cX}$ with $r > 0$, in which case we have by \eqref{FanoMirrorOfAnticanonical}
\begin{equation*}
\Theta_{k(- r K_{\cX})}(c_1(- r K_{\cX})) = r \Theta_{k(- r K_{\cX})}(c_1(- K_{\cX})) = r \cW(k \cA ) + O(k^{-1}).
\end{equation*}
Given this, we may perform our construction with respect to the polarisation $\cL = \cA$. Then, by definition, we have
\begin{align*}
[\theta_k] = k^{-1}\Theta^{-1}_{k \cA }(k c_1(\cA))|_{z=0} =  \Theta^{-1}_{k \cA }(c_1(\cA))|_{z=0} = r \cW(k \cA ) + O(k^{-1}).
\end{align*}
Thus, each term appearing in our expression \eqref{BaseLocusTerm} for the base locus contribution $\cB(\cW_k)$,
\begin{equation*}
\langle \res_{P_{n+1}}(\omega^{(k)}_+)\big| \res_{P_{n+1}}(\nabla_{\cW(k\cL)})^{-1}\big| \res_{P_{n+1}}(\omega^{(k)}_-)\rangle,\,\cD(P_{n+1}) \in \cD \setminus \cD_{\T}
\end{equation*}
is proportional to the quantity
\begin{align*}
&\res_{P_{n+1}}(\omega^{(k)}_+) = \res_{P_{n+1}}\left((\cW_k)^{-n} [\theta_k]^n \Omega_0\left(1 + O(k^{-1})\right)\right) \\
& = r^n \res_{P_{n+1}}\left( \Omega_0\left(1 + O(k^{-1})\right)\right)= O(k^{-1}),
\end{align*}
as required.
\end{proof}
\begin{exm}\label{NormalConeAnticanExm} Consider the degeneration to the normal cone $(\cX, k\cL_r) \cong (\cX, -\frac{k}{2} K_{\cX})$ of a point in $\PP^1$ with parameter $r = \frac{1}{2}$. We have
\begin{equation*}
\cW_k = \cW(k\cL_r) = x + \frac{e^{-2\pi k}}{x} + \frac{e^{-2\pi k}}{x'} + x' + e^{\pi k} x x',
\end{equation*}
The mirror map does not involve quantum corrections, giving
\begin{align*}
\theta_k = \frac{1}{2} \cW_k,\, \psi_k =  x + \frac{e^{-2\pi k}}{x} +  e^{\pi k} x x'
\end{align*} 
(see \cite{Stoppa_LargeComplex}, Example 2.9). As we saw, $\cW_k$, $\theta_k$ and $\psi_k$ extend to explicit anticanonical pencils on the toric del Pezzo $S_6 \subset \PP[x_0: \cdots: x_6]$ given by the toric boundary $S_6 \cap \{x_0 = 0\}$ together with the anticanonical sections
\begin{align*}
&\cW_k = x_1 + e^{-2\pi k} x_4 +  e^{-2\pi k} x_6 + x_3 + e^{ \pi k} x_2,\,\theta_k =  \frac{1}{2} \cW_k,\\
&\psi_k =  x_1 +  e^{-2\pi k} x_4 +  e^{ \pi k} x_2    
\end{align*} 
(we use the same notation for the sections and the corresponding pencils). 

The toric boundary of $S_6$ is given by $C_i = \PP[x_i, x_{i+1}]$ (with $i = 1, \ldots, 6$ cyclic) and so the torus fixed points are given by $p_i := C_i \cap C_{i+1}$. Thus, the only torus fixed point of $S_6$ contained in the base locus of $\cW_k$ is $p_4 = C_4 \cap C_5$. There are torus fixed points $p'_4 := T_{p_4} C_4$, $p''_4 := T_{p_4} C_5$ of $\barcT_{\T} \to S_6$ mapping to $p_4$. The values of $\tilde{\theta}$, $\tilde{\psi}$ at $p'_4$, $p''_4$ can be computed by taking the limit along $[x_4: 1] \in C_4$ as $x_4 \to 0$, respectively $[1: x_6] \in C_5$ as $x_6 \to 0$. We compute
\begin{align*}
&\tilde{\theta}|_{p} \equiv \frac{1}{2},\,\tilde{\psi}|_{p_i} =  1,\, p_i = p_1, p_3, p'_{4}, p_6,\,\tilde{\psi}|_{p_i} =  0,\,p_i = p_2, p''_{4}, p_5,\\ 
&\tilde{\theta}\left(\frac{n c}{n+1} \tilde{\theta} - 1 + \tilde{\psi}\right)|_{p} = \frac{1}{2}\left(- \frac{1}{2} + \tilde{\psi}\right)|_{p}. 
\end{align*}

Note that the divisors along which $\tilde{\psi} \equiv 1$ are given by the proper transforms of $C_1$ and $C_4$, which are contained in $\cX^{\vee, 1}$, $\cX^{\vee, 2}$ respectively.

It remains to compute the other factors appearing in the residue theorem. For this, near a fixed point $p_{i-1} = C_{i-1} \cap C_{i}$ for $i\neq 5$, we can cover a neighbourhood of a boundary component $C_{i} \subset S$ with local holomorphic coordinate patches 
\begin{equation*}
(z_{(i)}:= \frac{x_0}{x_{i}}, w_{(i)} = \frac{x_{i+1}}{x_{i}}),
\end{equation*}
such that $C_i$ is cut out by $z_{(i)} = 0$. For instance, near $p_1 \in C_2$, in coordinates $(z_{(2)}, w_{(2)})$, we have
\begin{align*}
&\Omega_0 = d\log x \wedge d\log x' = d\log \frac{x_1}{x_0} \wedge d\log \frac{x_3}{x_0}\\
&= d\log \frac{1}{w_{(2)}} \wedge d\log \frac{w_{(2)}}{z_{(2)}} = - d\log z_{(2)} \wedge d\log w_{(2)},
\end{align*}
where we used the toric relation $x_1 x_3 = x_0 x_2$. Similarly,  
\begin{align*}
&\cW_k = \frac{x_1}{x_0} + e^{-2\pi k} \frac{x_4}{x_0} +  e^{-2\pi k} \frac{x_6}{x_0} + \frac{x_3}{x_0} + e^{ \pi  k} \frac{x_2}{x_0}\\
& = \frac{1}{w_{(2)}} + e^{-2\pi k} w_{(2)}  +  e^{-2\pi k} \frac{z_{(2)}}{w_{(2)}} + \frac{w_{(2)}}{z_{(2)}} + e^{ \pi   k} \frac{1}{z_{(2)}},
\end{align*}
using the relations $x_1 x_4 = x^2_0$, $x_3 x_6 = x_0^2$, from which we compute
\begin{equation*}
\del_{z_{(2)}} \log \cW_k = \frac{-w_{(2)} e^{3 \pi  k }-e^{2 \pi  k} w^2_{(2)}+z^2_{(2)}}{z_{(2)} \left(w_{(2)} e^{3 \pi  k }+e^{2 \pi 
   k} \left(w^2_{(2)}+z_{(2)}\right)+z \left(w^2_{(2)}+z_{(2)}\right)\right)}.
\end{equation*}
It follows that we have $\Omega_0|_{p_1} = -1,\,\res(\nabla_{\cW_k})\big|_{p_1} = -1$.

Similar computations show that in fact we have $\Omega_0|_{p_i} = -1,\,\res(\nabla_{\cW_k})\big|_{p_i} = -1$ at all torus fixed points (including $p'_4$, $p''_4$).

According to \eqref{DFDecomposition}, we have
\begin{equation*}
\DF\big(\cX, -\frac{1}{2} K_{\cX}\big) = \sum_p\frac{1}{2}\left(- \frac{1}{2} + \tilde{\psi}\right)|_{p} + \cB(\cW_k) + O(k^{-1}) = \frac{1}{4} + \cB(\cW_k) + O(k^{-1}).
\end{equation*}
On the other hand, by the intersection-theoretic formula, 
\begin{align*}
&\DF\big(\cX, -\frac{1}{2} K_{\cX}\big) = \frac{1}{2}\DF\big(\cX, - K_{\cX}\big) = \frac{1}{2}\left(\frac{1}{2} (- K_{\cX})^2 + (K_{\cX} - \olo_{\PP^1}(-2)).(- K_{\cX})\right)\\
&= \frac{1}{2}\left(\frac{7}{2} - 7 +4\right) = \frac{1}{4}.
\end{align*}
This shows that we must have $\cB(\cW_k) = O(k^{-1})$, as required by Theorem \ref{VanishingThm}. 
\end{exm}
\begin{exm}\label{ProductAnticanExm} Consider the product test configuration  $(\cX, k\cL_r) \cong (\PP(\olo \oplus \olo(1)), -\frac{k}{3} K_{\cX})$ with parameter $r = \frac{1}{3}$. We have
\begin{equation*}
\cW_k = \cW(k\cL_r) = \frac{e^{-2\pi k}}{x x'} + x + x' + e^{\frac{2}{3}\pi k } x x'. 
\end{equation*}
The mirror map does not involve quantum corrections, giving
\begin{align*}
\theta_k = \frac{1}{3} \cW_k,\, \psi_k = x + x'. 
\end{align*} 

We know that $\cW_k$, $\theta_k$ and $\psi_k$ extend to explicit anticanonical pencils on the toric del Pezzo $S_6 \subset \PP[x_0: \cdots: x_6]$ given by the toric boundary $S_6 \cap \{x_0 = 0\}$ together with the anticanonical sections
\begin{align*}
&\cW_k =  e^{-2\pi k} x_5 + x_1 + x_3 + e^{\frac{2}{3}\pi k } x_2,\,\theta_k =  \frac{1}{3} \cW_k,\,\psi_k =  x_1 + x_3.    
\end{align*} 
The only torus fixed points of $S_6$ contained in the base locus of $\cW_k$ are $p_3 = C_3 \cap C_4$ and $p_5 = C_5 \cap C_6$. There are torus fixed points $p'_3 := T_{p_3} C_3$, $p''_3 := T_{p_3} C_4$ of $\barcT_{\T} = \Bl_{p_3, p_5} S_6 \to S_6$ mapping to $p_3$. The values of $\tilde{\theta}$, $\tilde{\psi}$ at $p'_3$, $p''_3$ can be computed by taking the limit along $[x_3: 1] \in C_3$ as $x_3 \to 0$, respectively $[1: x_5] \in C_4$ as $x_5 \to 0$. Similarly, there are torus fixed points $p'_5 := T_{p_5} C_5$, $p''_5 := T_{p_5} C_6$ of $\barcT_{\T}$ mapping to $p_5$, and the values of $\tilde{\theta}$, $\tilde{\psi}$ at $p'_5$, $p''_5$ can be computed by taking the limit along $[x_5: 1] \in C_5$ as $x_5 \to 0$, respectively $[x_1 : 1] \in C_6$ as $x_1 \to 0$.

We compute
\begin{align*}
&\tilde{\theta}|_{p} \equiv \frac{1}{3},\,\tilde{\psi}|_{p_i} =  1,\, p_i = p_2, p'_3, p''_5, p_6,\,\tilde{\psi}|_{p_i} =  0,\,p_i = p_1, p''_3, p_4, p'_5,\\ 
&\tilde{\theta}\left(\frac{n c}{n+1} \tilde{\theta} - 1 + \tilde{\psi}\right)|_{p} = \frac{1}{3}\left(\frac{1}{2} \frac{(\frac{1}{3}c_1(\PP^1)).(c_1(\PP^1))}{(\frac{1}{3}c_1(\PP^1))^2} \frac{1}{3} - 1 + \tilde{\psi}\right)|_{p}= \frac{1}{3}\left( - \frac{1}{2} + \tilde{\psi}\right)|_{p}. 
\end{align*}

The divisors along which $\tilde{\psi} \equiv 1$ are the proper transforms of $C_3$, $C_6$, contained in the central fibres of $\cX^{\vee, 1}$, $\cX^{\vee, 2}$ respectively.

Computing in local coordinates as in the previous example we also find $\Omega_0|_{p_i} = -1,\,\res(\nabla_{\cW_k})\big|_{p_i} = -1$ at all torus fixed points (including $p'_3$, $p''_3$, $p'_5$, $p''_5$).

According to \eqref{DFDecomposition}, we have
\begin{equation*}
\DF\big(\cX, -\frac{1}{3} K_{\cX}\big) = \sum_p\frac{1}{3}\left(- \frac{1}{2} + \tilde{\psi}\right)|_{p} + \cB(\cW_k) + O(k^{-1}) = \cB(\cW_k) + O(k^{-1}).
\end{equation*}
On the other hand, since $\cX$ is a product, we know that $\DF\big(\cX, -\frac{1}{3} K_{\cX}\big) = 0$. This shows that we must have $\cB(\cW_k) = O(k^{-1})$, as required by Theorem \ref{VanishingThm}. 
\end{exm}
\begin{exm}\label{SlopeUnstableSurfaceExmVanishing}Following Example \ref{SlopeUnstableSurfaceExm} for $X = \Bl_p \PP^2$, we have \emph{singular} total spaces of test configurations $\cX^{\vee, i}$, $i = 1, \ldots, 4$, for the \emph{orbifold} $T_{P^{\circ}}$ ``mirror" to $X$ (i.e. the toric, singular intersection of two quadrics in $\PP^4$), interchanged by a natural symmetry of $\cW_{\cX}(-k K_{\cX})$ (compatible with all our constructions), $\cX^{\vee, 1} \cong \cX^{\vee, 2}$, $\cX^{\vee, 3} \cong \cX^{\vee, 4}$. Assuming that the stationary phase formula holds for $\cW_{\cX}(-k K_{\cX})$ on $\barcT_{\T} = \widetilde{T}_{(\cP')^{\circ}}$, Theorem \ref{VanishingThm} provides smooth twisted formal test configurations $(\widetilde{\cX}^{\vee, i}, [\eta^i_k], [\xi^i_k])$ (i.e. without positivity conditions on $[\eta^i_k]$, $[\xi^i_k]$) such that
\begin{align*}
\DF(\cX, -K_{\cX}) = \sum^m_{i=1}\DF( \widetilde{\cX}^{\vee, i}, [\eta^i_k], [\xi^i_k]) + O(k^{-1}).
\end{align*}
By choosing $(\widetilde{\cX}^{\vee, i}, [\eta^i_k], [\xi^i_k])$ appropriately and taking some limits in $H^{1, 1}(\widetilde{\cX}^{\vee, i})$, we expect that this equality also induces the identity 
\begin{align*}
0 > \DF(\cX, -K_{\cX}) = 2(\DF(\cX^{\vee, 1}, [\hat{\eta}^1_k], [\hat{\xi}^1_k]) + \DF(\cX^{\vee, 3}, [\hat{\eta}^3_k], [\hat{\xi}^3_k])) + O(k^{-1}),
\end{align*}
where $\DF(\cX^{\vee, j}, [\hat{\eta}^j_k], [\hat{\xi}^j_k])$, $j = 1, 2$ are defined by the usual intersection-theoretic formula on the normal varieties $\cX^{\vee, j}$. In particular, at least one of $(\cX^{\vee, j}, [\hat{\eta}^j_k], [\hat{\xi}^j_k])$ ``destabilises" the orbifold $T_{P^{\circ}}$. 
\end{exm}
\begin{proof}[Proof of Theorem \ref{MainThmVariant}] We proceed as in the proofs of Theorems \ref{MainThmIntro} and \ref{VanishingThm}, by replacing the test configurations $(\cX^{\vee, i}, v_i, \eta^i_k, \xi^i_k)$ with $(\cX^{\vee, 0}, v_0, \cZ^0_k, \tilde{\cZ}^{0}_k)$,  $(\cX^{\vee, \infty}, v_{\infty}, \cZ^{\infty}_k, \tilde{\cZ}^{\infty}_k)$ as in Definition \ref{TestConfigDef}, but where now $v_0$, $v_{\infty}$ are \emph{fixed} (generically) as in Proposition \ref{FibrationsProp} $(iii)$, while only $(\cZ^0_k, \tilde{\cZ}^{0}_k)$,  $(\cZ^{\infty}_k, \tilde{\cZ}^{\infty}_k)$ are to be determined. Similarly to the proof of Theorem \ref{MainThmIntro}, we introduce functions $H_{0}, H_{\infty}\!: Z(v_i) \to \C$, $K_{0}, K_{\infty}\!: Z(v_i) \to \C$, now for \emph{fixed} $v_0$, $v_{\infty}$, such that
\begin{align*}
& H_{0}|_{\widetilde{Z}(v_0)} := d(v_0) f\tilde{\theta},\,\,H_{0}|_{Z(v_0)\setminus\widetilde{Z}(v_0)}:= 0,\\
& K_{0}|_{\widetilde{Z}(v_0)} := - t(v_0) - d(v_0) f \left(1 - \tilde{\psi} \right),\,\,K_{0}|_{Z(v_0)\setminus\widetilde{Z}(v_0)}:= 0,   
\end{align*}
respectively
\begin{align*}
& H_{\infty}|_{\widetilde{Z}(v_{\infty})} := d(v_{\infty}) f\tilde{\theta},\,\,H_{\infty}|_{Z(v_{\infty})\setminus\widetilde{Z}(v_0)}:= 0,\\
& K_{\infty}|_{\widetilde{Z}(v_{\infty})} := - t(v_{\infty}) - d(v_{\infty}) f \left(1 - \tilde{\psi} \right),\,\,K_{\infty}|_{Z(v_{\infty})\setminus\widetilde{Z}(v_{\infty})}:= 0,   
\end{align*}
and we observe that the analogue of \eqref{PrescribedH}, namely
\begin{align*} 
&(\cZ^0_k)'(p) = H_0(p),\,\,(\tilde{\cZ}^0_k)'(p) = K_0(p), \textrm{ for all } p \in Z(v_0),\\
&(\cZ^{\infty}_k)'(p) = H_{\infty}(p),\,\,(\tilde{\cZ}^{\infty}_k)'(p) = K_{\infty}(p), \textrm{ for all } p \in Z(v_0),
\end{align*} 
can always be solved for equivariant representatives $(\cZ^0_k)', (\cZ^{\infty}_k)' \in H^*_{S^1}(\cX^{\vee}, \C)$ of classes $\cZ^0_k, \cZ^{\infty}_k \in H^*(\cX^{\vee}, \C)$, since the cohomology ring $H^*(\cX^{\vee}, \C)$ localises to the fixed point set under the $\C^*$-action. The rest of the arguments in the proofs of Theorems \ref{MainThmIntro}, \ref{VanishingThm} apply without changes.
\end{proof}

\addcontentsline{toc}{section}{References} 
 
\bibliographystyle{abbrv}
 \bibliography{biblio_Kstab}

\noindent SISSA, via Bonomea 265, 34136 Trieste, Italy;\\
Institute for Geometry and Physics (IGAP), via Beirut 2, 34151 Trieste, Italy\\
jstoppa@sissa.it    
\end{document}